\documentclass{article}
\usepackage{amsmath}
\usepackage{latexsym}
\usepackage{amssymb}
\usepackage{amsthm}
\usepackage{verbatim}
\usepackage{enumerate}
\font\logic=msam10 at 10pt
\newcommand{\forces}{\mbox{\logic\char'015}}
\newcommand{\restrict}{\mbox{\logic\char'026}}

\def\undertilde#1{{\baselineskip=0pt\vtop
  {\hbox{$#1$}\hbox{$\scriptscriptstyle\sim$}}}{}}

\newtheorem{thrm}{Theorem}[section]
\newtheorem{lem}[thrm]{Lemma}
\newtheorem{cor}[thrm]{Corollary}

\newtheoremstyle{hdefinition}%
  {\topsep}%
  {\topsep}%
  {\upshape}%    Definitions, etc not in italics.
  {}%
  {\bfseries}%
  {.}%           would be better with period.
  { }%
  {\thmnumber{#2 }\thmname{#1}\thmnote{ \rm(#3)}}%

\theoremstyle{hdefinition}
\newtheorem{remark}[thrm]{Remark}

\newtheorem{df}[thrm]{Definition}
\newtheorem{propty}[thrm]{Property}

\begin{document}

\title{Bounding by canonical functions, with CH\thanks{This
research was conducted while both authors were in residence
at the Mittag-Leffler Institute. We thank
the Institute for its hospitality.}}
\author{Paul Larson\thanks{The research of the first author was supported
in part by the NSERC grants of Juris Stepr\={a}ns, Paul Szeptycki
and Franklin D. Tall, and the Centre de Recerca Matem\`{a}tica of
the Institut d'Estudis Catalan.} \and Saharon Shelah\thanks{The
research of the second author was supported by the Israel Science
Foundation, founded by the Israel Academy of Sciences. Publication
number 746.}}

\pagenumbering{arabic}
\maketitle

\begin{abstract}
We show that the members of a certain class of
semi-proper iterations do not add countable sets
of ordinals. As a result,
starting from suitable large cardinals
one can obtain a model
in which the Continuum Hypothesis holds and every
function from $\omega_{1}$ to $\omega_{1}$ is bounded
on a club by a canonical function for an ordinal less
than $\omega_{2}$.
\end{abstract}

%\begin{abstract} We show the consistency of the continuum hypotheis
%with the statement that every function from $\omega_{1}$ to
%$\omega_{1}$ is bounded by a canonical function on a club.
%\end{abstract}

\section{Introduction}

Given an ordinal $\gamma$, a function $f \colon \omega_{1}
\to Ord$ is a \emph{canonical function for} $\gamma$ if the
empty condition (i.e., $\omega_{1}$) in the forcing
$\mathcal{P}(\omega_{1})/NS_{\omega_{1}}$ forces that $j(f)(\omega^{V}_{1}) =
\gamma$, where $j$ is the elementary embedding induced by the
generic. For each $\alpha < \omega_{1}$, the constant function with
value $\alpha$ is the canonical function for $\alpha$. For
$\alpha \in [\omega_{1}, \omega_{2})$, a canonical function $f$ for
$\alpha$ is obtained by taking a bijection $g \colon \omega_{1}
\to \alpha$ and letting $f(\beta)$ be the ordertype of
$g[\beta]$.
%The first interesting questions about canonical functions are whether
%there is a canonical function for $\omega_{2}$, and if there is,
%whether it is the constant function $\omega_{1}$.
In this paper we let \emph{Bounding} denote the statement that
every function from $\omega_{1}$ to $\omega_{1}$ is bounded on a club
subset of $\omega_{1}$ by a canonical function for an ordinal
less than $\omega_{2}$. It is fairly easy to see that if the
nonstationary
ideal on $\omega_{1}$ ($NS_{\omega_{1}})$ is saturated, then Bounding holds.
The second author has shown \cite{Shf} that given
the existence of a Woodin cardinal there is a semi-proper forcing
making $NS_{\omega_{1}}$ saturated, and it has been known for some time that
there is a simpler forcing making Bounding hold from a weaker
large cardinal hypothesis.
%In each case, one shows that the corresponding
%forcing is semi-proper by finding for each suitable elementary substructure
%a corresponding expanded structure with the same countable ordinals such that
%one can construct a condition generic for the expanded structure.
%One essential difference in the two cases is that in forcing saturation this
%expanded structure is found by contradiction, whereas in the second case
%it is constructed directly (see Theorem \ref{std}).
%This implies that $\omega_{1}$ is
%the $\omega_{2}$-nd canonical function.
The most quotable result in this paper is that this standard forcing to make
every function from $\omega_{1}$ to $\omega_{1}$ bounded by a canonical
function is $(\omega, \infty)$-distributive (i.e., it does not add
$\omega$-sequences of ordinals), and so
this statement is
consistent with the Continuum Hypothesis, even in the presence of
large cardinals. This is in contrast with
saturation, as Woodin \cite{W} has shown that if $NS_{\omega_{1}}$ is
saturated and sufficiently large cardinals exist, then there
is a definable counterexample to CH.
We give a more general theorem
stating that the members of a certain class of semi-proper iterations
are $(\omega, \infty)$-distributive. This class includes the standard
forcing
to make Bounding hold, and is general
enough to show that a generalization of Bounding for certain
sets of reals is also
consistent with CH, answering a question in \cite{W}.

The key construction used in the proof generalizes the notion of
$\alpha$-properness from
Chapter V of \cite{Shf} to semi-proper forcing. Briefly, a
forcing is $\alpha$-semi-proper if for any $\in$-chain of
countable elementary submodels of length
$\alpha$, there is a condition which
is  simultaneously semi-generic for each model in the
sequence. The problem in applying the
method to show that a given
improper iteration is $(\omega, \infty)$-distributive
is that for a given model $N$ in the
sequence,
%the base model
%$N_{0}$, for which we would like to find a condition
%$p$ forcing a value to
%$N_{0}[G] \cap \undertilde{G}$ (or
%$\undertilde{N}[G_{0}]\restrict (N_{0} \cap \omega_{1})$,
%in the case of improper forcing),
%is that
%$N[G_{p_{\alpha}}] \cap X \cap length(\bar{Q})$
$N[G] \cap \kappa$ can be a proper
superset of $N \cap \kappa$,
where $\kappa$ is the length of the iteration and $G$ is generic for
some initial segment,
so that new steps appear.
For the forcings in this paper, however,
we have a good understanding of how
to enlarge each such $N$, as well as how to
produce the appropriate tower of models to
overcome this.

This can be generalized further, getting the
consistency of certain forcing axioms, using ideas from
\cite{Shf}, Chapters V and VIII, and \cite{Sh656}. The reader is
referred to \cite{Sh311} for more on this topic and on RCS in
particular.

%Woodin \cite{W} has shown that if $NS_{\omega_{1}}$ is saturated and
%large cardinals exist, then CH fails. One motivation for
%the work in this paper was to see whether this theorem holds with
%the weaker hypothesis of bounding by canonical functions.
%We have shown that it doesn't.

Interest in this question derives also from the study of Woodin's
$\mathbb{P}_{max}$ forcing \cite{W}, which produces a model in
which CH fails and all forceable $\Pi_{2}$ sentences for
$H(\omega_{2})$ hold simultaneously. It is not known whether all
such $\Pi_{2}$ sentences forceably consistent with CH can hold
together with CH. The generalized form of bounding in this paper
is a candidate for showing that this is impossible. Candidates for
the other half of the incompatibility appear in \cite{L, Sh638}
and in Section 10.6 of \cite{W}.

%In the forcing order, we let $q \leq p$ mean than $q$ is more
%informative than $p$.

\section{Skolem Hulls}

Given a structure $M$ with a predicate $<^{*}$ for a wellordering of the
domain of $M$, and a subset $X$ of the domain of $M$, we let
$Sk_{(M, \in, <^{*})}(X)$ denote the Skolem hull of $X$ in $(M, \in, <^{*})$.
If $X$ is a countable elementary substructure of some
$H(\chi)$ and $\eta$ is an ordinal in $X$, then
we let $D^{X}_{\eta}$ be the set of all $a \in [\eta]^{<\omega}$
such that $f(a) \in X$ for all $f \colon [\eta]^{|a|} \to \omega_{1}$ in
$X$.
The point is that if $<^{*}$ is a wellordering of $H(\chi)$,
$X \prec (H(\chi), \in, <^{*})$
is countable, and  $z$ is a subset of
some ordinal $\eta \in X$, then
 $$Sk_{(H(\chi), \in, <^{*})}(X \cup z)
 = \{ f(a) : a \in [z]^{<\omega} \wedge f \in H(\chi)^{([\eta]^{<\omega})}
\cap X\}$$
and so if $z$ is such that $[z]^{<\omega} \subset
D^{X}_{\eta}$,
$$Sk_{(H(\chi), \in, <^{*})}(X \cup z)
\cap \omega_{1} = X \cap \omega_{1}.$$
Note that if $\eta < \eta'$ are ordinals in an elementary
submodel $X$, then $$D^{X}_{\eta} = D^{X}_{\eta'} \cap [\eta]^{<\omega}.$$

In order to verify that our forcings satisfy the semi-properness
condition we require, we need to repeatedly apply a certain
simultaneous-extendibility property. The following lemma will be
useful in this regard.

\begin{lem}\label{helper} Say that $M$ is a structure and $<^{*}$ is a
wellordering of the domain of $M$. Let $X$ be a countable elementary
substructure of $(M, \in, <^{*})$. Let $\eta$ be an ordinal in $X$ and let
$z_{0}, z_{1}$ be countable subsets of $\eta$ such that $[z_{0}]^{<\omega}
\subset D^{X}_{\eta}$ and $[z_{1}]^{<\omega} \subset D^{X'}_{\eta}$, where
$X' = Sk_{(M, \in, <^{*})}(X \cup z_{0})$. Then $[z_{0} \cup
z_{1}]^{<\omega} \subset D^{X}_{\eta}$.
\end{lem}

\begin{proof} Let $a \in [z_{0}]^{<\omega}$ and $b\in [z_{1}]^{<\omega}$,
and let $f \colon [\eta]^{|a| + |b|} \to \omega_{1}$. Then $f_{a}
\colon [\eta]^{|b|} \to \omega_{1}$, defined by letting $f_{a}(y)
= f(a\cup y)$, is a function in
$$X' = Sk_{(M, \in, <^{*})}(X \cup z_{0}),$$
so $f(a\cup b) \in X' \cap \omega_{1} = X \cap \omega_{1}$.
\end{proof}

Say that $X,Y$ are countable elementary submodels of some
$(H(\chi), \in, \leq_{\chi})$ with $X \in Y$. Let $\eta$ be an
ordinal in $X$ and let $\gamma \in Y \cap \eta$. Then even though
$Sk_{(H(\chi), \in, \leq_{\chi})}(X \cup \{ \gamma \})$ is not
directly definable in $Y$, the set itself is in $Y$, since it is
equivalent to $\{ f(\gamma) \mid f \in X \wedge f \colon \eta \to
H(\chi) \}$.

\begin{remark}\label{keytoo} Similarly, many of the arguments in this paper prove
facts about sequences of elementary submodels by induction on the
length of the sequence. Of course if $X,Y$ are countable
elementary submodels of $(H(\chi), \in, \leq_{\chi})$ with $X \in
Y$, $Y$ does not see that $X$ is an elementary submodel, and so
the induction hypothesis cannot be applied directly in $Y$.
However, the statement about whether an object exists with a
certain relation to $X$, a semi-generic extending a certain
condition, say, is formalizable in $Y$, and so if one exists in
$H(\chi)$ then one exists in $Y$.
\end{remark}

\section{A class of $(\omega, \infty)$-distributive iterations}

Each step of the iterations we are considering is a forcing
which shoots a continuous increasing sequence
of length $\omega_{1}$ through a
given stationary set of countable sets of ordinals. Under certain
assumptions on the stationary set, such forcings are
a typical example of improper forcings which preserve stationary
subsets of $\omega_{1}$, and they have been well studied in recent
years (see, for example, \cite{FeJ}). For
the iterations in this paper, we require
that these stationary sets be definable from sets from the ground model
with the help of functions from $\omega_{1}$ to $\omega_{1}$
added by initial segments of the iteration.

\begin{df}\label{aitdef}  Let $\langle \lambda_{\rho} : \rho < \kappa
\rangle$ be a continuous increasing sequence of ordinals,
and let  $\mathcal{A} =
\langle A^{\rho}_{\beta} :
\rho < \kappa, \beta < \omega_{1} \rangle$ be such that each
$A^{\rho}_{\beta} \subset [\lambda_{\rho+1}]^{<\omega_{1}}$.
Given $\rho < \kappa$ and $f \colon \omega_{1}\to\omega_{1}$, let
$Q_{\rho, f}$ be the forcing whose conditions are countable, continuous,
increasing sequences $\langle x_{\beta} : \beta \leq \gamma \rangle$
such that for each $\beta \leq \gamma$, $x_{\beta} \cap \omega_{1} \in
\omega_{1}$ and $x_{\beta} \in A^{\rho}_{f(x_{\beta} \cap \omega_{1})}$,
ordered by extension.
Then an $\mathcal{A}$-\emph{iteration} is a structure
$\bar{Q} =
\langle P_{\rho}, \undertilde{Q}_{\rho}, \undertilde{f}_{\rho} : \rho <
\kappa\rangle$ such that
\begin{enumerate}
\item[(a)] $\langle P_{\rho}, \undertilde{Q}_{\rho}
 : \rho < \kappa\rangle$ is a Revised Countable Support iteration,
%\item[(b)] each $P_{\rho}$ has cardinality (or, formally, density) $\leq
%\lambda_{\rho}$ for nonlimit $\rho$,  and $\leq 2^{\lambda_{\rho}}$ for
%$\rho$ a successor.
\item[(b)] each $\undertilde{f}_{\rho}$ is a $P_{\rho}$-name for a
function from $\omega_{1}$ to $\omega_{1}$,
%for $\beta = 0$ or
%$\beta < \rho$,
%\item[(c)] if $\undertilde{f}$ is a $P_{\rho}$-name for a function from
%$\omega_{1}$ to $\omega_{1}$, then for some $\beta \in [\rho, \kappa)$,
%$\undertilde{f}_{\beta} = \undertilde{f}$,
\item[(c)] for all $\rho < \kappa$, $1_{P_{\rho}}
\forces \undertilde{Q}_{\rho} = Q_{\rho, \undertilde{f}_{\rho}}$.
\end{enumerate}
\end{df}

The stationary sets in our iterations must also satisfy certain
extendibility conditions with respect to the countable elementary submodels
of
a sufficiently large initial segment of the universe.

\begin{thrm}\label{both} Let $\langle \lambda_{\rho} : \rho < \kappa
\rangle$ be a continuous increasing sequence of strong limit
cardinals with supremum $\kappa$. Fix a regular cardinal $\chi >
(2^{\kappa})^{+}$,  and let $\leq_{\chi}$ be a wellordering of
$H(\chi)$. Let $\mathcal{A} = \langle A^{\rho}_{\beta} : \rho <
\kappa, \beta < \omega_{1} \rangle$ be such that each
$A^{\rho}_{\beta} \subset [\lambda_{\rho+1}]^{<\omega_{1}}$, and
let $\bar{Q} = \langle P_{\rho}, \undertilde{Q}_{\rho},
\undertilde{f}_{\rho} : \rho < \kappa\rangle$ be an
$\mathcal{A}$-iteration, such that the following hold.
\begin{enumerate}
\item\label{up} For all $\rho < \kappa, \beta < \omega_{1}$,
if $E \in
A^{\rho}_{\beta}$ and $E' \in [\lambda_{\rho+1}]^{<\omega_{1}}$
with $E \subset E'$
%and $E \cap \omega_{1} = E' \cap \omega_{1}$
then $E' \in A^{\rho}_{\beta}$, \item\label{upto} For all $\rho <
\kappa$, $1_{P_{\rho}}$ forces that for all countable $X \prec
(H(\chi)^{V^{P_{\rho}}}, \in, \leq_{\chi})$ with
$\lambda_{\rho+1}, \bar{Q} \in X$, and for all $\beta <
\omega_{1}$ there exists a countable $z \subset \lambda_{\rho+1}$
such that $[z]^{<\omega}$ is a subset of
$$\bigcap
\{ D^{Z}_{\lambda_{\rho+1}} : Z \prec (H(\chi)^{V^{P_{\rho}}},
\in, \leq_{\chi}) \wedge \lambda_{\rho+1}, \bar{Q} \in Z \in X\}$$
and, letting $Y = Sk_{(H(\chi)^{V^{P_{\rho}}}, \in,
\leq_{\chi})}(X \cup z)$,
\begin{enumerate}
%\item[(i)] $X \subset Y$
\item[(i)] $X \cap \omega_{1} = Y \cap
\omega_{1}$,
\item[(ii)] $Y \cap \lambda_{\rho+1} \in A^{\rho}_{\beta}$.
\end{enumerate}
\item\label{upthre} For all countable $X \prec (H(\chi), \in,
\leq_{\chi})$ with $\bar{Q} \in X$, and for all $\rho \in X \cap
\kappa, \beta < \omega_{1}$ there exists a countable $Y \prec
(H(\chi), \in, \leq_{\chi})$ such that
\begin{enumerate}
\item[(i)] $X \subset Y$
\item[(ii)]\label{beet} $X \cap V_{\lambda_{\rho} + 2} = Y \cap
V_{\lambda_{\rho} + 2}$,
\item[(iii)] $Y \cap \lambda_{\rho+1} \in A^{\rho}_{\beta}$.
\end{enumerate}
\end{enumerate}
Let $P$ be the Revised Countable Support limit of
$\langle P_{\rho}, \undertilde{Q}_{\rho} : \rho < \kappa\rangle$.
%corresponding to $\bar{Q}$
Then $P$ is $(\omega, \infty)$-distributive.
\end{thrm}

\begin{remark} By standard RCS arguments, if in the statement of
Theorem \ref{both} we assume
in addition that
$\kappa$ is strongly  inaccessible, then $P$ is $\kappa$-c.c. (see
Theorem \ref{dumbtwo}).
\end{remark}

Conditions \ref{upto} and \ref{upthre} of Theorem \ref{both}
could easily be subsumed into one condition; indeed, in our
applications we verify both conditions at the same time.
The conditions correspond
to separate parts of the proof, however. Condition \ref{upto} is
needed to show that each successor step of the iteration is
$\alpha$-semi-proper for all countable $\alpha$
(Theorem \ref{yugotprops}), and could in fact
be replaced by this requirement, though this would mean more
work in applying the theorem. Condition \ref{upthre} is used to
construct the systems of models which we use to show that
no countable sets of ordinals are added (Lemma \ref{issys}).

%Since each $Q_{\rho, f}$ is semi-proper, by the semi-properness
%preservation theorems in \cite{Shf}, the forcing in Theorem
%\ref{both} is semi-proper, if one uses Revised Countable Support
%instead of Countable Support.
%Indeed, iterations as in Definition \ref{aitdef} are typically improper,
%and so when semi-proper are usually iterated using Revised Countable
%Support. In our case, the iterations add no countable sets of
%ordinals, so the two
%approaches give rise to identical forcings. The Countable Support
%approach requires two additional steps. In section \ref{csi},
%we discuss Countable Support limits of $(\omega, \infty)$-distributive
%forcings. These theorems are special cases
%of the RCS theory in \cite{Sh311, Shf}. Furthermore, the
%Countable Support approach requires a simultaneous induction
%showing that the iteration is $(\omega, \infty)$-distributive and
%semi-proper, instead of, as in the RCS approach, asserting that
%the entire iteration is semi-proper, and then showing
%$(\omega, \infty)$-distributivity by induction. Otherwise the
%two approaches are the same.

Since in the end we show that $P$ does not add $\omega$-sequences,
these are actually Countable Support iterations. Nonetheless, at
this time we do not have a proof of (the corresponding version of)
Theorem \ref{both} which avoids RCS.

%and so the whole proof can be
%rearranged without mention of RCS.
%We have defined it this way for expositional convenience.
%Carrying out the argument with Countable Support requires simultaneous
%induction on the preservation of $\alpha$-semi-properness and not
%adding $\omega$-sequences, instead of,
%as we have done it here, asserting that $\alpha$-semi-properness is
%preserved for the whole iteration and then proving by induction
%that $\omega$-sequences aren't added. Otherwise, the two approaches are
%the same.

%One could strengthen the theorem by relaxing the upwards closure of
%the sets $A^{\rho}_{\beta}$, requiring that Condition \ref{up}
%above holds only for some $\eta < \lambda_{\rho}$, instead of for
%$\omega_{1}$. This would also require modifying the definition of
%$D^{X}_{\eta}$. Notice that Condition \ref{beet} implies that each
%$\lambda_{\rho+ 1} > 2^{\lambda_{\rho}}$. There are several other
%technical improvements one could make, though at this point we
%see no application for them.

For the rest of this paper, sets denoted by $\kappa$, $\langle
\lambda_{\rho} : \rho < \kappa \rangle$, $\chi$, $\leq_{\chi}$,
$Q_{\rho, f}$, $\bar{Q}$, $P$, $\mathcal{A}$ and $\langle
A^{\rho}_{\beta} : \rho < \kappa, \beta < \omega_{1} \rangle$ are
supposed to have the properties given in the hypotheses of Theorem
\ref{both}. This policy will be modified  in two ways. First, in
the application sections, we will add  extra properties for these
terms. On the other hand, in proving the main  theorem we
sometimes state our lemmas more generally in terms of certain
properties of these objects, temporarily forgetting the others.
We hope that it is clear when we are doing this, so that there
will be no confusion.

\section{Proof of Theorem \ref{both}}

The proof of Theorem \ref{both} breaks into two largely disjoint
parts, which come together only in Lemma \ref{alldone}. In the
first part, we define $\alpha$-semi-proper forcing, where $\alpha$
is a countable ordinal, and show that for each countable ordinal
$\alpha$, $\alpha$-semi-properness is preserved under RCS
iterations. In Theorem \ref{yugotprops} we use Condition
\ref{upto} from Theorem \ref{both} to show that each successor
step from one of these iterations is forced to be
$\alpha$-semi-proper for all countable $\alpha$, which then
carries over to the entire iteration. In the second part of the
proof we use the end-extension property given in Condition
\ref{upthre} of Theorem \ref{both} to find a system of elementary
submodels suitable for proving $(\omega, \infty)$-distributivity,
which is shown in Lemmas \ref{mgensex}, \ref{fixone}, \ref{issys}
and \ref{alldone}. By Theorem \ref{otwo}, the iteration is
$\kappa$-c.c.

%We wish to
%see that the iteration $P$ is $\kappa$-c.c. and
%$(\omega, \infty)$-distributive.

%To save ourselves some
%extra cases in Lemmas \ref{dumbone}, \ref{cspres} and
%\ref{alldone}, we will prove these by simultaneous
%induction on $\kappa$, though those lemmas should also hold if
%the condition of $(\omega, \infty)$-distributivity is
%removed.

%Theorem \ref{dumbone} then
%gives the induction step for semi-properness, and,
%letting $P_{\kappa} = RCS(P)$,

\subsection{Revised Countable Support}\label{csi}

The reader is referred to \cite{Shf, Sh311} for the
definition of Revised Countable Support and its basic analysis.
Alternate presentations of RCS can be found in
\cite{DoFu, Mi}. We will show that the arguments here
follow in all three versions.
The only facts about RCS that we
need outside of this section are Theorems
\ref{boxone}, \ref{otwo} and \ref{datsit},
and Property \ref{propone}.
These are (essentially - see Remark \ref{miproprem})
proved in \cite{Mi, Mi2} for
the version of RCS presented in \cite{Mi}. We
present proofs of these facts that work for the
versions of RCS in \cite{Shf} and \cite{DoFu}.

The following two theorems hold for all
presentations of RCS.
The first gives an important property of RCS which
distinguishes it from Countable Support, that
names for conditions are essentially conditions
themselves.

\begin{thrm}\label{boxone} Say that $\langle P_{\alpha},
\undertilde{Q}_{\alpha} : \alpha < \kappa \rangle$ is an RCS
iteration with RCS limit $P$. Fix $p \in P$ and $\gamma < \kappa$.
Suppose that $A$ is a maximal antichain in $P_{\gamma}$ below $p
\restrict \gamma$ and $f : A \to P$ is a function such that for
each $a \in A$ $f(a) \leq p$ and $f(a) \restrict \gamma = a$. Then
there is a condition $p' \leq p$ such that $p' \leq  p$, $p
\restrict \gamma = p' \restrict \gamma$ and each $a \in A$ forces
that $p' \restrict [\gamma, \kappa) = f(a) \restrict [\gamma,
\kappa)$.
\end{thrm}

Theorem \ref{dumbone} follows from Theorem \ref{boxone}
and Corollary 2.8 from Chapter X of \cite{Shf} (and the
corresponding theorem from \cite{DoFu}).

\begin{thrm}\label{dumbone} Suppose
that $\langle P_{\gamma},
\undertilde{Q}_{\gamma} : \gamma < \kappa \rangle$ is an
RCS iteration with RCS limit $P_{\kappa}$ such that
each  $\undertilde{Q}_{\gamma}$ is forced to be
semi-proper and
%\item each $P_{\gamma}$ is $(\omega, \infty)$-distributive,
each $\undertilde{Q}_{\gamma}$ forces the corresponding
$P_{\gamma}$ to have cardinality $\aleph_{1}$.
Then the following hold.
\begin{enumerate}
\item $P_{\kappa}$ is semi-proper.
%\item for all $\gamma \leq \delta \leq \kappa$,
%if $p \in X \cap P_{\kappa}$ and
%$q \in P_{\gamma}$ is $(X, P_{\gamma})$-semi-generic with $q \leq p \restrict
%\gamma$, then there exists an $(X, P_{\delta})$-semi-generic condition $q' \in
%P_{\delta}$ below $p\restrict \delta$ such that $q' \restrict \gamma = q$
\item\label{oboyoboy} For all regular cardinals $\chi > 2^{|P_{\kappa}|}$,
for all countable $X
\prec H(\chi)$
with $P_{\kappa} \in X$, for all $\gamma < \delta$ in $X \cap (\kappa+1)$, if
\begin{enumerate}
\item $p \in P_{\kappa}$,
\item $q$ is an $(X, P_{\gamma})$-semi-generic
condition in $P_{\gamma}$ below $p \restrict \gamma$,
\item $q$ forces ``there exists $r$ in $P_{\kappa} \cap X$
such that $p \restrict [\gamma, \kappa) =  r \restrict [\gamma,
\kappa)$,"
\end{enumerate}
then there exists an $(X, P_{\delta})$-semi-generic
condition $q' \in P_{\delta}$
such that
$q' \restrict \gamma = q$ and $q' \leq p \restrict \delta$.
\end{enumerate}
\end{thrm}

In \cite{DoFu}, RCS is characterized by the following
property.

\begin{propty}\label{propone} If $\langle P_{\alpha} :
\alpha < \kappa \rangle$ is an RCS iteration with
RCS limit $P$
and $p \in P$, then for
all $q \leq p$ in $P$ there exist
$\gamma < \kappa$, $r \in P_{\gamma}$ such
that $r \leq q \restrict \gamma$ and either
$r \forces cof(\check{\kappa}) = \omega$ or
$r \forces \forall \alpha \in (\check{\gamma},
\check{\kappa})\text{ }p(\alpha) = 1_{\undertilde{Q}_{\alpha}}$.
\end{propty}

If $\langle P_{\alpha} :
\alpha < \kappa \rangle$ is an RCS iteration with
RCS limit $P$ then for each condition $p$ in $P$
there is an associated $P$-name $supp(p)$ for the
support of $p$, the set of $\alpha < \kappa$
such that $$p(\alpha) \neq 1_{\undertilde{Q}_{\alpha}},$$
as decided by the $P$-generic filter. One can
easily prove by induction on $\kappa$ that
Property \ref{propone} implies that for
each $p$ in $P$ $supp(p)$ is forced
to be countable.

\begin{remark}\label{miproprem} Property \ref{propone} is
shown in \cite{Mi} for the version of RCS in \cite{Mi}, for
the special case where $\kappa = \omega_{1}$ (so the
first possibility for $r$ cannot hold), en route
to (essentially) proving Theorem \ref{dumbtwo}. Together
these two facts give Property \ref{propone} for this
version of RCS, since in the remaining case
(where $|P_{\alpha}| < cof(\kappa)$ for all $\alpha <
\kappa$) a maximal antichain deciding the supremum
of $supp(p)$ for some condition $p$ must have
maximal restriction in some initial segement of
the iteration.
\end{remark}

In \cite{Shf}, we have the following property
(see Definition
1.1 and Claim 1.3 (1) of \cite{Shf}, Chapter X).

\begin{propty}\label{proptwo} If $\langle P_{\alpha} :
\alpha < \kappa \rangle$ is an RCS iteration of limit length with
RCS limit $P$ and $p \in P$, then there exist $P$-names $\tau_{i}$
$(i < \omega)$ for ordinals less than $\kappa$ such that $$1_{P}
\forces supp(\check{p}) = \{ \tau_{i} : i < \omega \}$$ and such
that for each $i < \omega$ and each $q \leq p$ in $P$ there exist
$\gamma < \kappa, r \in P_{\gamma+1}$ such that $r \leq q
\restrict (\gamma+1)$ and $r \forces \tau_{i} = \check{\gamma}$.
\end{propty}

If we assume that in addition each $\undertilde{Q}_{\alpha}$
forces that the cardinality of the corresponding $P_{\alpha}$ is
$\aleph_{1}$, then Property \ref{proptwo} implies Property
\ref{propone}. We need only check the limit case. Fix $p$ and $q$,
let $\{ \tau_{i} : i < \omega \}$ be the sequence given by
Property \ref{proptwo} with respect to $p$, and assume that there
is no pair $\gamma, r$ such that
\begin{itemize}
\item $\gamma < \kappa$,
\item $r \in P_{\gamma}$,
\item $r \leq q \restrict \gamma$,
\item $r \forces cof(\check{\kappa}) = \omega$.
\end{itemize}
If some $P_{\alpha}$ has cardinality greater than the cofinality
of $\kappa$, then (assuming that $q$ is in the generic filter),
$\kappa$ will have cofinality $\omega_{1}$ in the extension by
$P_{\alpha+1}$. Then, working in the $P_{\alpha+1}$-extension, let
$X$ be a countable elementary substructure of a large enough
$H(\chi)$ with $P, p,q, \{\tau_{i} : i < \omega\} \in X$. Define
$\zeta = \bigcup(X \cap \kappa)$ and let $q' \leq q \restrict
[\alpha + 1, \kappa)$ be $(X, P/P_{\alpha+1})$-semi-generic. Then
$q'$ forces that $supp(p)$ will be contained in $\zeta$. This
means that for each $i < \omega$ every condition $q'' \leq q'$ in
$P/P_{\alpha + 1}$ will be compatible with an $r \in
P_{\zeta}/P_{\alpha+1}$ forcing that $\tau_{i} < \zeta$. But then
$q' \restrict \zeta$ will also have this property, and thus $q'
\restrict \zeta$ also forces that $supp(p) \subset \zeta$.
Therefore, there exists a condition $r \in P_{\gamma}$ for some
$\gamma < \kappa$ such that $r\restrict (\alpha + 1)$ forces that
$\gamma$ will have the properties of $\zeta$ as above, and that $r
\restrict [\alpha+1, \gamma)$ will satisfy the properties of $q'
\restrict \zeta$. Such an $r$ suffices.

On the other hand, if $cof(k) > |P_{\alpha}|$ for
all $\alpha$, then there is a sequence of
pairs $(\gamma_{i}, r_{i})$ $(i < \omega)$ such
that
\begin{itemize}
\item the $\gamma_{i}$'s are increasing,
\item each  $r_{i} \in P_{\gamma_{i}}$ forces a bound
below $\kappa$ on $\tau_{i}$,
\item for all $i < j$, $r_{j} \restrict \gamma_{i} = r_{i}$,
\item each $r_{i} \leq p \restrict \gamma_{i}$.
\end{itemize}
Then the limit of the $r_{i}$'s is the desired condition.
To find $\gamma_{i+1}, r_{i+1}$, find a maximimal
antichain $A$ in $P_{\gamma_{i}}$ below $r_{i}$ such
that for each $a \in A$ there is a
$\gamma_{a} < \kappa$ and an $r_{a} \in P_{\gamma_{a}}$
such that
\begin{itemize}
\item $r_{a} \restrict \gamma_{i} = a$, \item $r_{a} \leq p$,
\item $r_{a} \forces \tau_{i} = \check{\gamma}_{a}$.
\end{itemize}
Then apply Theorem \ref{boxone} to $A$ and the
function $a \mapsto r_{a}$ to find $r_{i+1}$,
and let $\gamma_{i+1} = sup\{ \gamma_{a} :
a \in A\}$, which must be below $\kappa$ since
$cof(\kappa) > |P_{\gamma_{i}}|$.

Similar considerations give the following (see also
Lemma 36.5 of \cite{J}). The point is that the
RCS limit of an iteration of cofinality $\omega_{1}$
is just the direct limit, and the ordinals of
cofinality $\omega_{1}$ are stationary below
$cof(\kappa)$ as below. Then for any maximal
antichain $A$ in $P$ as below
there is some $\gamma < \kappa$ of
cofinality $\omega_{1}$ such that
$A \cap P_{\gamma}$ is a maximal antichain in
$P_{\gamma}$. But then $|A| \leq |P_{\gamma}|$.

\begin{thrm}\label{dumbtwo} Say that $\langle P_{\gamma},
\undertilde{Q}_{\gamma} : \gamma < \kappa \rangle$ is an RCS iteration with RCS limit
$P$ such that $cof(\kappa) > |P_{\gamma}|$ for
all $\gamma < \kappa$, and such that each
$\undertilde{Q}_{\gamma}$ forces the corresponding $P_{\gamma}$ to
have cardinality $\aleph_{1}$.  Then $P$ is $cof(\kappa)$-c.c.
\end{thrm}

%For the rest of this section,
%we fix the items supposed in the statement of Theorem \ref{both}.
%Nonetheless, we sometimes state our lemmas more generally in terms of certain
%properties of $\chi, \lambda_{\rho}, \kappa$, forgetting temporarily
%that we have fixed them. It should be clear when we are doing this, and when
%we are referring to our fixed values.

%Condition \ref{upto} in
%Theorem \ref{both} gives that each of the forcings
%$\undertilde{Q}_{\rho}$ is semi-proper over the
%$P_{\rho}$-extension.
Applying Theorems \ref{dumbone} and \ref{dumbtwo},
%the preservation theorems in \cite{Shf}
we have reduced the proof of Theorem \ref{both} to showing
that each $P_{\beta}$ $(\beta < \kappa)$ is $(\omega, \infty)$-distributive.

\begin{thrm}\label{otwo} Let $\langle P_{\beta}, \undertilde{Q}_{\beta} :
\beta < \kappa \rangle$ be an RCS iteration of strongly
inaccessible length with RCS limit $P_{\kappa}$ such that
\begin{itemize}
\item each $1_{P_{\beta}}$ forces the corresponding
$\undertilde{Q}_{\beta}$ to
be semi-proper,
\item each $\undertilde{Q}_{\beta}$ makes the corresponding
$P_{\beta}$ have cardinality $\aleph_{1}$,
\item each $|P_{\beta}| < \kappa$,
\item each $P_{\beta}$ is $(\omega, \infty)$-distributive.
\end{itemize}
Then for all $\rho < \beta \leq \kappa$, $P_{\beta}/P_{\rho}$ is semi-proper
and $\kappa$-c.c. Therefore $P_{\kappa}$ is $(\omega, \infty)$-distributive.
\end{thrm}

%It is shown in \cite{Sh311} that if $\kappa$ is strongly inaccessible
%then an RCS limit of length $\kappa$ of forcings of cardinality less
%than $\kappa$ is $\kappa$-c.c.. The rest of Claim \ref{otwo} follows from
%\cite{Shf}, Chapter X.

\subsection{Semi-generics for sequences}

The following are generalizations of
ideas from Chapters V, X and XII of
\cite{Shf}.

\begin{df} Let $\alpha$ be a countable ordinal.
\begin{enumerate}
\item The set $SEQ_{\alpha}(\chi)$ consists of all $\bar{N}
= \langle N_{\beta} : \beta < \alpha \rangle$ such that
%(in \cite{Shf},V, called...)
\begin{enumerate}
%\item $\rho < \omega_{1}$,
%\item $\bar{N} = \langle N_{\beta} : \beta < \alpha \rangle$,
\item each $N_{\beta}$ is a countable elementary substructure of
$(H(\chi), \in, \leq_{\chi})$,
%\item $\beta < \gamma < \alpha \Rightarrow N_{\beta} \prec N_{\gamma}$,
\item for each $\gamma < \alpha$, $\bar{N}\restrict \gamma =
\langle N_{\beta} : \beta < \gamma \rangle \in N_{\gamma}.$
%\item $\beta < \rho \Rightarrow \bar{N}\restrict (\beta + 1) \in
%N_{\beta + 1}.$ SHELAH'S ORIGINAL
\end{enumerate}
\item Let $\bar{N} \in SEQ_{\alpha}(\chi)$ and let $P$ be a forcing
construction in $N_{0}$. A condition $p
\in P$ is $(\bar{N}, P)$-\emph{semi-generic} if
$$p \forces \check{N}_{\delta}[\undertilde{G}_{P}] \cap \omega_{1} =
\check{N}_{\delta} \cap \omega_{1}$$
for all $\delta < \alpha$, where $\undertilde{G}_{P}$ is a
$P$-name for the generic filter.
The condition $p$ is $(\bar{N} \restrict [\gamma, \beta],
P)$-\emph{semi-generic} if this holds
for all $\delta \in [\gamma, \beta]$.
\item A forcing construction $P$ is $\alpha$-\emph{semi-proper} if for
every $\bar{N} \in SEQ_{\alpha}(\chi)$ with $P \in N_{0}$ and for all $p
\in P \cap N_{0}$ there is
an $(\bar{N}, P)$-semi generic $q \in P$ such that $p \geq q$.
%and recall (same).
\end{enumerate}
\end{df}

Given $\bar{N} \in SEQ_{\alpha}(\chi)$ and $P \in N_{0}$, if
$G \subset P$ is a generic filter then $\bar{N}[G]$ is the
sequence $\langle N_{\beta}[G] : \beta < \alpha \rangle$.
Both parts of the following lemma are immediate.  Note that
in the second part there is no need to find a $p \in \undertilde{G}_{P}$
which is $(\bar{N}, P)$-semi-generic.

\begin{lem}\label{helpme}
Let $\alpha$ be a countable ordinal and fix $\bar{N} \in
SEQ_{\alpha}(\chi)$.
\begin{enumerate}
\item\label{cooper} Let $P$ be an $\alpha$-semi-proper forcing in
$N_{0}$, and
let $\undertilde{Q}$ be a $P$-name in $N_{0}$
for an $\alpha$-semi-proper forcing.
Let $q \in P$  be $(\bar{N}, P)$-semi-generic and let
$\undertilde{p}$ be a $P$-name for a condition in $(\bigcup
\bar{N}[\undertilde{G}_{P}]) \cap \undertilde{Q}$. Then there is
a $P$-name $\undertilde{r}$ for a condition in $\undertilde{Q}$
such that $1_{P} \forces \undertilde{r} \leq \undertilde{p}$ and
$(q,\undertilde{r})$ is $(\bar{N}, P*Q)$-semi-generic.
\item $1 \forces_{P} \check{\bar{N}}[\undertilde{G}_{P}] \in
SEQ_{\check{\alpha}}(\check{\chi})^{V[\undertilde{G}_{P}]}.$
\end{enumerate}
\end{lem}

%The proof of the theorem below is a straightforward combination of the proofs
%in \cite{Shf} that CS iterations preserve $\rho$-properness and
%RCS iterations preserve semi-properness.

\begin{lem}\label{cspres} Let $\alpha$ be a countable ordinal.
Suppose that $P_{\kappa}$ is the RCS limit of an RCS
iteration $\langle P_{\rho}, \undertilde{Q}_{\rho} : \rho <
\kappa \rangle$ in a set forcing extension $V^{*}$
such that each $1_{P_{\rho}}$ forces the corresponding
$\undertilde{Q}_{\rho}$ to be $\alpha$-semi-proper, and
%each $P_{\rho}$ $(\rho < \kappa)$ is
%$(\omega, \infty)$-distributive and
each $\undertilde{Q}_{\rho}$ forces $P_{\rho}$ to have
cardinality $\aleph_{1}$.
Then the following hold in $V^{*}$.
\begin{enumerate}
\item $P_{\kappa}$ is $\alpha$-semi-proper.
\item\label{tweetie} Fix $\beta \leq \gamma \leq \kappa$,
and fix $\bar{N} \in SEQ_{\alpha}(\chi)$ with
$P_{\kappa}$ and $\beta, \gamma \in \bigcup\bar{N}$.
Let $p \in P_{\kappa}$, $q \in P_{\beta}$ and
$\delta < \alpha$ be such that
\begin{itemize}
\item $\beta \in N_{\delta}$,
\item $p$ is $(\bar{N}\restrict \delta, P_{\kappa})$-semi-generic,
\item $q$ is $(\bar{N}\restrict [\delta, \alpha),
P_{\beta})$-semi-generic,
\item $p \restrict \beta \geq q$,
\item $q$ forces that for some $r \in P_{\kappa} \cap N_{\delta}$,
$p \restrict [\beta, \kappa) = r \restrict
[\beta, \kappa)$.
\end{itemize}
Let $\delta^{*} \in [ \delta, \alpha)$ be such that $\gamma \in N_{\delta^{*}}$.
Then there exist $q^{*} \in
P_{\gamma}$ and $p^{*} \in P_{\kappa}$  such that
\begin{itemize}
\item $p^{*}$ is $(\bar{N}\restrict \delta^{*}, P_{\kappa})$-semi-generic,
\item $q^{*}$ is $(\bar{N}\restrict[\delta^{*}, \alpha),
P_{\gamma})$-semi-generic,
\item $q^{*} \leq p^{*} \restrict \gamma$,
\item $p^{*} \leq p$,
\item $q^{*} \restrict \beta = q$,
\item $q^{*}$ forces that for some $r \in P_{\kappa} \cap N_{\delta^{*}}$,
$p^{*} \restrict [\gamma, \kappa) = r \restrict
[\gamma, \kappa)$.
\end{itemize}
\end{enumerate}
\end{lem}

%\begin{lem}\label{rcspres} Let $\alpha$ be a countable ordinal.
%If $\bar{P} =
%\langle P_{\beta}, \undertilde{Q}_{\beta} : \beta < \gamma
%\rangle$ is an RCS limit of $\alpha$-semi-proper forcings such that
%for each $\beta < \gamma$ there exists $n \in \omega$ such that
%$$1_{P_{\beta + n}} \forces |P_{\beta}| \leq \aleph_{1},$$ then
%$\bar{P}$ is $\alpha$-semi-proper.
%\end{lem}

\begin{proof} We prove part
\ref{tweetie} by induction on $\alpha$, working in $V^{*}$.
The first part follows immediately. We fix the notation that
$G_{\rho}$ is the generic filter for $P_{\rho}$, for $\rho <
\kappa$. The case $\alpha = 1$ is
given by Theorem \ref{dumbone}. For the case where
$\alpha = \alpha' + 1$, there are two subcases. If
$\delta^{*} < \alpha'$, then we may assume by the induction
hypothesis that $\delta^{*} = \delta$ and that there is
a $(\bar{N}\restrict [\delta^{*}, \alpha'), P_{\gamma})$-semi-generic
$q' \in P_{\gamma}$ such that $q' \restrict \beta = q$ and
$q' \leq p \restrict \gamma$. Then since $q$ forces that
$N_{\alpha'}[G_{\beta}]$
will be elementary in $H(\chi)^{V^{*}[G_{\beta}]}$, where
$G_{\beta}$ is the generic filter for $P_{\beta}$, we
can replace $q'$ with a $q''$ with the additional property
that $q$ forces that $q''\restrict [\beta, \gamma)$ will be
equal to $r \restrict [\beta, \gamma)$ for some $r \in N_{\alpha'}
\cap P_{\kappa}$.
By part \ref{oboyoboy} of Theorem \ref{dumbone},
then, there is a $q^{*} \in P_{\gamma}$ such that
$q^{*} \restrict \beta = q$, $q^{*} \leq q''$ and $q^{*}$ is
$(N_{\alpha'}, P_{\gamma})$-semi-generic. Such a $q^{*}$
suffices.

For the subcase $\delta^{*} = \alpha'$, by the induction
hypothesis
there is a $(\bar{N}\restrict \alpha', P_{\kappa})$-semi-generic
$q' \in P_{\gamma}$ such that $q' \restrict \beta = q$ and
$q' \leq p$. Then since $q$ forces that
$N_{\alpha'}[G_{\beta}]$
will be elementary in $H(\chi)^{V^{*}[G_{\beta}]}$, where
$G_{\beta}$ is the generic filter for $P_{\beta}$, we
can replace $q'$ with a $q''$ with the additional property
that $q$ forces that $q''\restrict [\beta, \kappa)$ will be
equal to $r \restrict [\beta, \kappa)$ for some $r \in N_{\alpha'}
\cap P_{\kappa}$, and we can let $p^{*} = q''$.
Then as in the previous paragraph, by
part \ref{oboyoboy} of Theorem \ref{dumbone}
there is a $q^{*}$ as desired.

The case where $\alpha$ is a limit and $\kappa = \kappa'+1$ is similar,
now fixing $\alpha$ and inducting on $\kappa$.
By the induction hypothesis for $\alpha$ we may assume that
$\delta = \delta^{*}$
%there is a $(\bar{N}\restrict \delta^{*}, P_{\kappa})$-semi-generic
%$p$ such that $p' \leq p$, and
and by the induction hypothesis for $\kappa$ we may assume that
there is a $(\bar{N}\restrict [\delta^{*}, \alpha), P_{\kappa'})$-semi-generic
$q' \in P_{\kappa'}$ such that $q' \restrict \beta = q$ and
$q' \leq p \restrict \kappa'$. Since $\undertilde{Q}_{\kappa'}$ is
forced to be $\alpha$-semi-proper,
%and since
%$q'$ forces that $N_{\delta^{*}}[G_{\kappa'}]$
%will be elementary in $H(\chi)^{V^{*}[G_{\kappa'}]}$, where
%$G_{\kappa'}$ is the generic filter for $P_{\kappa'}$,
by Theorem \ref{boxone} (or part \ref{cooper} of
Lemma \ref{helpme}) there is a $(\bar{N},
P_{\kappa})$-semi-generic
$q^{*}$ below $p$ such
that $q^{*} \restrict \kappa = q'$.

For the case where $\alpha$ and $\kappa$ are both limits there are
two subcases (in each of which we may assume that $\delta =
\delta^{*}$). If $\kappa$ has cofinality $\omega$ or $\omega_{1}$,
or if $|P_{\rho}| < cof(\kappa)$ for each $\rho < \kappa$, then we
fix an increasing sequence $\langle \xi_{i} : i < \omega \rangle$
cofinal in $\alpha$, with $\xi_{0} = \delta$. If the cofinality of
$\kappa$ is countable, fix an increasing sequence of ordinals
$\langle \beta_{i} : i < \omega \rangle \in N_{0}$ cofinal in
$\kappa$, with $\beta_{0} = \beta$. Otherwise, let $\beta_{i}$ be
any ordinal in $N_{\xi_{i}}$ greater than $sup(\bigcup \bar{N}
\restrict \xi_{i} \cap \kappa)$, again with $\beta_{0} = \beta$.
Let $p_{0}=p$ and let $q_{0} = q$. Alternately choose conditions
$p_{i+1}, q_{i+1}$ $(i < \omega)$, such that
\begin{enumerate}
\item each $p_{i+1}$ is a condition in $P_{\kappa}$,
%\item each $p_{i+1}\restrict \beta_{i} = p_{i} \restrict \beta_{i}$,
\item each $q_{i+1}$ is a condition in $P_{\beta_{i+1}}$,
\item each $q_{i} \leq p_{i} \restrict \beta_{i}$,
\item each $p_{i+1} \leq p_{i}$,
%\item each $q_{i}$ forces $p_{i+1}\restrict [\beta_{i}, \kappa) \leq
%p_{i} \restrict [\beta_{i}, \kappa)$,
\item for all $i < j < \omega$, $q_{j} \restrict \beta_{i} = q_{i}$,
\item each $p_{i+1}$ is $(\bar{N}\restrict\xi_{i+1}, P_{\kappa})$-semi-generic,
%\item in the case that $cof(\kappa)$ is countable, each $q_{i+1}$ is $(\bar{N},
%P_{\beta_{i}})$-semi-generic,
\item each $q_{i+1}$ is
$(\bar{N} \restrict [\xi_{i+1}, \alpha), P_{\beta_{i+1}})$-semi-generic,
%(here $(\bar{N}, P_{\beta_{i+1}})$-semi-generic doesn't make sense,
%since $\beta_{i+1} \not\in N_{0}$)
\item each $q_{i}$ forces that for some condition $r \in
P_{\kappa} \cap N_{\xi_{i}}$, $r \restrict [\beta_{i}, \kappa) =
p_{i} \restrict [\beta_{i}, \kappa)$.
\end{enumerate}
For the case where no condition in any $P_{\eta}$ makes
$cof(\kappa) \leq \omega_{1}$, we modify condition 8
as follows:
\begin{enumerate}
\item[8a.] each $q_{i}$ forces that for some  condition $r \in
P_{\kappa} \cap N_{\xi_{i}}$ such that for some $\gamma < \kappa$
$1_{P_{\kappa}}\forces supp(\check{r}) \subset \check{\gamma}$, $r
\restrict [\beta_{i}, \kappa) = p_{i} [\beta_{i}, \kappa)$.
\end{enumerate}

That such conditions exist is immediate by the induction
hypothesis (8a follows from Property \ref{propone}). Then the
limit of the $q_{i}$'s (call it $q^{*}$) will be the desired
$(\bar{N}, P_{\kappa})$-semi-generic, as long as it is below each
$p_{i}$. This fact follows from the fact that $q^{*}$ forces that
$\{ \beta_{i} : i < \omega \}$ will be cofinal in $\bigcup \{
supp(p_{i}) : i < \omega \}$. This is clear if $cof(\kappa) =
\omega$, and if $cof(\kappa) = \omega_{1}$ it follows from the
fact that each $N_{\xi_{i}} \cap \kappa$ will be cofinal in
$N_{\xi_{i}}[G_{\beta_{i}}] \cap \kappa$ since $q_{i}$ is
$(N_{\xi_{i}}, P_{\beta_{i}})$-semi-generic. For the remaining
case it follows from condition 8a.

If $|P_{\rho}| \geq cof(\kappa)$ for some $\rho < \kappa$, then
we may assume that $\beta > \rho$ and so $cof(\kappa) \leq
\omega_{1}$ in the $P_{\beta}$-extension. Then we may apply
the previous argument in the $P_{\beta}$-extension, along
with Theorem \ref{boxone}.
\end{proof}

\begin{thrm}\label{datsit} Let $\alpha$ be a countable ordinal.
Say that $P_{\kappa}$ is the RCS limit of an RCS
iteration $\langle P_{\beta},
\undertilde{Q}_{\beta}
: \beta < \kappa \rangle$ such
%that each $P_{\beta}$ is $(\omega, \infty)$-distributive,
that each $P_{\beta+1} = P_{\beta} * \undertilde{Q}_{\beta}$
forces $|P_{\beta}| \leq \aleph_{1}$ and
each $1_{P_{\beta}}$ forces that $\undertilde{Q}_{\beta}$ is
$\alpha$-semi-proper. Then for all $\beta < \gamma \leq
\kappa$, $P_{\gamma}/P_{\beta}$ is
$\alpha$-semi-proper for every $\alpha < \omega_{1}$.
\end{thrm}

To show that each $Q_{\rho, f}$ is $\alpha$-semi-proper for
all $\alpha < \omega_{1}$, we show that we
can extend sequences of models in a suitable way.

\begin{lem}\label{yugotnun} Assume that
$\lambda < \chi, \bar{Q}$, and
$\langle A_{\beta} :  \beta < \omega_{1} \rangle \subset
[\lambda]^{<\omega_{1}}$
are such that
\begin{itemize}
\item for all $\beta < \omega_{1}, E \in A_{\beta}$, if $E \subset
E'$ then $E' \in A_{\beta}$, \item for all countable $X \prec
(H(\chi), \in, \leq_{\chi})$ with $\lambda, \bar{Q} \in X$, and
for all $\beta < \omega_{1}$ there exists a countable $z \subset
\lambda$ such that $[z]^{<\omega}$ is a subset of
$$\bigcap
\{ D^{Z}_{\lambda} : Z \prec (H(\chi), \in, \leq_{\chi}) \wedge
\bar{Q}, \lambda \in Z \in X\}$$ and, letting $Y = Sk_{(H(\chi),
\in, \leq_{\chi})}(X \cup z)$,
\begin{enumerate}
%\item[(i)] $X \subset Y$
\item[(i)] $X \cap \omega_{1} = Y \cap
\omega_{1}$,
\item[(ii)] $Y \cap \lambda \in A_{\beta}$.
\end{enumerate}
\end{itemize}
Let $\bar{N} \in SEQ_{\xi}(\chi)$ (for some $\xi < \omega_{1}$)
with $\bar{Q}, \lambda \in N_{0}$, and let
%\item each member of $\bar{N}$ is in $C$,
$g\colon\omega_{1} \to \omega_{1}$ be a function such
that for each $\eta < \xi$ $g \restrict (\cup\{ N_{\eta'} \cap \omega_{1} :
\eta' < \eta\} + 1) \in N_{\eta}$.
Then
%(actually, $\gamma
%= \lambda$ is o.k. (except $|N''| > \aleph_{0}$)
there exists
$\bar{N}' \in SEQ_{\xi}(\chi)$
such that
\begin{enumerate}
%\item each $N'_{\eta} \cap \lambda_{\rho + 1}$
%end-extends $N_{\eta} \cap \lambda_{\rho + 1}$,
\item[(a)] each $N_{\eta} \subset N'_{\eta}$,
\item[(b)] each $N'_{\eta} \cap \omega_{1} = N_{\eta} \cap
\omega_{1}$,
\item[(c)] each $N'_{\eta} \cap \lambda \in
A_{g(N'_{\eta} \cap \omega_{1})}$,
\item[(d)]\label{cat} letting $N' = \bigcup\{ N'_{\eta} : \eta < \xi\}$,
$N' \cap \lambda \in A_{g(N' \cap \omega_{1})}$.
\end{enumerate}
\end{lem}

\begin{proof} By induction on $\xi$.
If $\bar{N}$
has a last model, then  we are done by
the extension assumption on countable elementary submodels of
$H(\chi)$.
For the limit case,
applying the induction hypothesis we may assume
that every condition already holds, except part (d).
Let $$N = \bigcup\{ N_{\eta} : \eta < \xi\}.$$
Applying the end-extension assumption, let
$z$ be a countable subset of $\lambda$ such that
$$[z]^{<\omega} \subset \bigcap
\{ D^{Z}_{\lambda} : Z \prec (H(\chi), \in, \leq_{\chi}) \wedge
\lambda, \bar{Q} \in Z \in N\}$$ and such that, letting
$N'=Sk_{(H(\chi), \in, \leq_{\chi})}(N \cup z)$, we have
\begin{enumerate}
\item $N \subset N'$,
\item $N \cap \omega_{1} = N' \cap \omega_{1}$,
\item $N' \cap \lambda \in
A_{g(N' \cap \omega_{1})}$.
\end{enumerate}
Let $h \colon\omega \to z$ be a bijection, and fix an increasing
sequence $\langle \zeta_{i} : i < \omega \rangle$ cofinal in $\xi$
with $\zeta_{0} = 0$. Now let each $$N'_{\eta} = Sk_{(H(\chi),
\in, \leq_{\chi})}(N_{\eta} \cup h[i_{\eta}]),$$ where $i_{\eta}$
is the largest $i$ such that $\zeta_{i} \leq \eta$. Since these
are all finite extensions, each initial sequence of the new
sequence is an element of the later models. By the choice of $z$,
each $N_{\eta} \cap \omega_{1} = N'_{\eta} \cap \omega_{1}$, and
so $\bar{N}'$ is as desired.
\end{proof}

\begin{thrm}\label{yugotprops} If $P_{\rho}$ is
as in Theorem \ref{both} and
$f \colon \omega_{1} \to \omega_{1}$ is a function added by
$P_{\rho}$, then $Q_{\rho, f}$ is
$\alpha$-semi-proper in the extension by $P_{\rho}$
for all $\alpha < \omega_{1}$.
%and the $\rho$-semi-generic can be taken to be a condition $p$ such that
%each $p \cap N_{\eta} \in N_{\eta + 1}$.
\end{thrm}

\begin{proof} By induction on $\alpha$,
using the Lemma \ref{yugotnun}
and working in the extension by $P_{\rho}$.
Fix $\bar{N} \in
SEQ_{\alpha}(\chi)$ and $p \in Q_{\rho, f}$
with $p, \bar{Q}, \lambda \in N_{0}$,
and let $\bar{N}'$ be as in Lemma \ref{yugotnun}, with respect to $\langle
A^{\rho}_{\beta} :
\beta < \omega_{1} \rangle$. If $\bar{N}'$
has a final model (i.e., if $\alpha$ is a successor), then we can choose a
$(\bar{N}'\restrict (\alpha - 1), Q_{\rho, f})$-semi-generic $p'\in
N'_{\alpha-1}$ extending $p$ by the induction hypothesis (see Remark
\ref{keytoo}). Then since $$N'_{\alpha -1} \cap \lambda_{\rho+1} \in
A^{\rho}_{f(N'_{\alpha -1} \cap \omega_{1})},$$ any $N'_{\alpha
-1}$-generic for $Q_{\rho, f}$ extending $p'$ will suffice as the desired
condition. For the case where $\alpha$ is a limit, we let $\bar{N}$ and
$\bar{N}'$ be as in the proof of Lemma \ref{yugotnun}, and fix an
increasing sequence $\xi_{i}$ $(i <
\omega)$ of ordinals cofinal in $\alpha$.
Then letting
$p$ be $p_{0}$ we can successively pick conditions $p_{i}$ $(i < \omega)$
such that each $p_{i+1}$ is a $(\bar{N}'\restrict \xi_{i+1}, Q_{\rho,
f})$-semi-generic in $N'_{\xi_{i+1}}$ extending $p_{i}$, and such that the
last member of each $p_{i+1}$ contains $N'_{\xi_{i}} \cap \lambda_{\rho
+1}$. By Condition (d) from Lemma \ref{yugotnun}, the limit of the
$p_{i}$'s (adjoined by their union) will be a condition in $Q_{\rho, f}$,
and it will be $(\bar{N}, Q_{\rho, f})$-semi-generic since it extends
an $(\bar{N} \restrict \beta, Q_{\rho, f})$-semi-generic for each $\beta <
\alpha$.
\end{proof}

\subsection{Systems of models and conditions}

To show that each initial segment $P_{\beta}$ $(\beta <
\kappa)$
is $(\omega, \infty)$-distributive, we choose a suitable pair $(M, \bar{N})$
where $M$ is a countable elementary submodel of $H(\chi)$ and
$\bar{N} \in SEQ_{o.t.(M \cap \beta)}(\chi)$, and
find an $(\bar{N}, P_{\beta})$-semi-generic
condition in $P_{\beta}$ which extends an $M$-generic filter.
As defined below, a $(\bar{Q}, \rho, \beta)$-system is a partial
construction of such an object (minus the $M$-genericity requirement),
where $\rho$ is the length of the initial segment for the iteration
for which the desired condition has been constructed, and $\beta$ is the
target length.

\begin{df}\label{sysdef} Given ordinals $\rho \leq \beta \leq \kappa$,
we say that $(M, \bar{N}, p, q)$ is a $(\bar{Q},
\rho, \beta)$-\emph{system} if (letting $\delta = o.t.(M \cap \rho)$ and
$\gamma = o.t.(M \cap \beta)$):
\begin{enumerate}
\item $M \prec (H(\chi), \in, \leq_{\chi})$, \item $\bar{Q}$,
$\rho$, $\beta$ belong to $M$,
%\item $\rho \leq \beta$,
\item $\bar{N} \in SEQ_{\gamma}(\chi)$ and $\bar{Q} \in N_{0}$,
\item\label{sysfive} $q \in P_{\rho}$ and $q$ is
$(\bar{N}\restrict[\delta, \gamma), P_{\rho})$-semi-generic,
\item\label{syssix} for all $\rho' \in \rho \cap M$, $q \restrict \rho'$ is
$(N_{o.t.(M \cap \rho')}, P_{\rho'})$-semi-generic,
%\item $q$ forces a value to $\undertilde{G}_{p_{\rho}} \cap M$,
\item\label{g} for all $\eta \in M \cap \beta$,
\begin{enumerate}
\item\label{syssevena} for all $\nu \in N_{o.t(M \cap \eta)} \cap
\omega_{1}$, $M \cap \lambda_{\eta + 1} \in A^{\eta}_{\nu}$, \item
$M \cap V_{\lambda_{\eta}+2} \in  N_{o.t.(M \cap \eta)}$,
\end{enumerate}
\item\label{syslast} $p \in P_{\beta}$, $p \restrict \rho \geq q$
and for all $\eta \in [\rho, \beta) \cap M$, $p \restrict \eta \in
N_{o.t.(M \cap \eta)}$.
\end{enumerate}
\end{df}

\begin{remark}\label{sdog} Note that if $\beta = \rho + 1$ and
$(M, \bar{N}, p, q)$ is a $(\bar{Q}, \rho, \beta)$-system, then
$(M, \bar{N}, p, q)$ is also a $(\bar{Q}, \beta, \beta)$-system,
as (the second part of) Condition \ref{sysfive} becomes vacuous in the
second case, and Condition \ref{sysfive} in the first case completes
Condition \ref{syssix} in the second.
\end{remark}

%\begin{remark} Notice that if  $q \in P_{\rho}$ and $q$ is
%$(\bar{N}\restrict[\delta, \gamma), P_{\rho})$-semi-generic,
%then for each $\rho' \in M \cap (\rho + 1)$,
%$q\restrict \rho' \in P_{\rho'}$ is
%$(\bar{N}\restrict[\delta', \gamma), P_{\rho'})$-semi-generic,
%where $\delta' = o.t.(M \cap \rho')$.
%\end{remark}

\begin{df} Let $M \prec (H(\chi), \in, \leq_{\chi})$ and $\beta < \kappa$
with $P_{\beta} \in M$. A filter $g \subset P_{\beta} \cap M$ is $M$-\emph{generic} if
for all dense sets $D \subset P_{\beta}$ in $M$ the intersection
of $g$ and $D$ is nonempty. A condition $p \in M$ is a
\emph{potential}-$M$-\emph{generic} for
$P_{\beta}$ if
there exists a $g \subset P_{\beta} \cap M$ $M$-generic for $P_{\beta}$ such that
for all $\rho \in M \cap \beta$,
$1_{P_{\rho}}$ forces that if $$M \cap \lambda_{\rho + 1} \in
A^{\rho}_{\undertilde{f}_{\rho}(M \cap \omega_{1})}$$ then
$p(\rho) = g(\rho)^{\frown}\langle \cup g(\rho) \rangle,$ and if for
all $\rho \not\in M \cap \beta, p(\rho) = 1_{\undertilde{Q}_{\rho}}$.
\end{df}

Note that if $g \subset P_{\beta}$ is $M$-generic and $p$ is in
$g$ then since $supp(p)$ is a $P_{\beta}$-name for a countable
subset of $\beta$ the support of $p$
as determined by $g$ is contained
in $M \cap \beta$.
We need to see that potential-$M$-generics exist in suitable generality.

\begin{lem}\label{mgensex} For any $\beta \leq \kappa$, if
$(M, \bar{N}, p, \emptyset)$ is a
$(\bar{Q}, 0, \beta)$-system  with  $p \in M$,
then there is  a potential-$M$-generic $\bar{p}\leq p$ such that $\bar{p}
\restrict \eta \in N_{o.t.(\eta \cap M)}$ for all $\eta \in M \cap \beta$.
\end{lem}

\begin{proof}
Given $g \subset P \cap M$ $M$-generic for $P_{\beta}$, one
can easily build a corresponding potential-$M$-generic $p$, letting
$p(\rho)$ be the empty condition when
$$M \cap \lambda_{\rho + 1} \not\in
A^{\rho}_{\undertilde{f}_{\rho}(M \cap \omega_{1})}.$$
The point then is just to find an $M$-generic filter
$g \subset M \cap P_{\beta}$ with $p \in g$, such
that $g \cap P_{\eta} \in N_{o.t.(\eta \cap M)}$
for all $\eta \in  M \cap \beta$.

Inducting primarily on $\beta \in M \cap \kappa$, and secondarily
on $\gamma \in M \cap \beta$, we show that if $p^{*} \in P_{\beta}
\cap M$ and $g^{*} \in N_{o.t.(M \cap \gamma)}$ is an $M$-generic
filter for $P_{\gamma}$ such that $p^{*}\restrict \gamma \in
g^{*}$ and $g^{*} \cap P_{\gamma'} \in N_{o.t.(M \cap \gamma')}$
for all $\gamma' \in M \cap \gamma$, then there exists an
$M$-generic $g \subset P_{\beta}$ such that $g \cap P_{\gamma} =
g^{*}$, $p^{*} \in g$ and $g \cap  P_{\eta} \in N_{o.t.(M \cap
\eta)}$ for all $\eta \in M \cap \beta$. Note that for each
$\gamma \in M \cap \beta$, since
$$\mathcal{P}(\mathcal{P}(\lambda_{\gamma}))
\cap M \in N_{o.t.(M \cap \gamma)}$$
(and $|P_{\gamma}| \leq 2^{\lambda_{\gamma}}$),
the $M$-genericity of $g^{*}$ can be verified in
$N_{o.t.(M \cap \gamma)}$ ($M$ itself is not
in any member of $\bar{N}$, so this is not automatic).

%and that $g^{*} \subset P_{\gamma}$ in an $M$-generic
%filter in $N_{o.t.(M \cap \gamma)}$
%with $p^{*} \restrict \gamma \in g^{*}$ and $g^{*} \cap P_{\eta}
%\in N_{o.t.(M \cap \eta)}$ for all $\eta < \gamma$.
%Then $N_{o.t.(M \cap \gamma)}$ is the last member of $\bar{N}$. Since
%$$\bar{N} \restrict o.t.(M \cap \gamma) \in N_{o.t.(M \cap \gamma)}$$
%and $N_{o.t.(M \cap \gamma)}$ is an elementary submodel of $H(\chi)$,
%there is an $M$-generic filter $g^{*} \subset M \cap P_{\beta - 1}$
%in $N_{o.t.(M \cap \gamma)}$ as required, as computed in $N_{o.t.(M \cap
%\gamma)}$.

The base case $\beta = 0$ is trivial.
For the successor step, we may
assume that $\gamma + 1 = \beta$
Then since for the successor case there is no restriction on the
last coordinate of the $M$-generic filter (other than $M$-genericity),
$g^{*}$ can be extended in any fashion to an $M$-generic for
$P_{\beta}$.

For the limit, fix an increasing sequence
$\eta_{i}$ $(i < \omega)$ cofinal in $\beta \cap M$, and let $p = p_{0}$.
Now in $\omega$ many steps alternately pick
\begin{itemize}
\item $M$-generic $g_{i}\subset P_{\eta_{i}}$
in $N_{o.t.(M \cap \eta_{i})}$
(applying the induction hypothesis plus the elementarity of
$N_{o.t.(M \cap \eta_{i})}$
plus the fact that $\bar{N} \restrict o.t.(M \cap \eta_{i}) \in
N_{o.t.(M \cap \eta_{i})}$, as in Remark \ref{keytoo})
with $p_{i} \restrict
\eta_{i} \in g_{i}$ such that
\begin{itemize}
\item $g_{i} \cap P_{\eta} \in N_{o.t.(M \cap \eta)}$ for all
$\eta \in M \cap \eta_{i}$,
%\item satisfies the lemma for the case $\beta = \eta_{i}$ as
%computed in the model $N_{o.t.(M \cap \eta_{i})}$, and
\item $g_{i} \cap P_{\eta_{j}} = g_{j}$ for all $j < i$,
\end{itemize}
\item $p_{i+1} \leq p_{i}$ in $M$ meeting the $i$th dense set in $M$ for
$P_{\beta}$ such that
$p_{i+1} \restrict \eta_{i} \in g_{i}$ (such a $p_{i+1}$
exists because $g_{i}$ is $M$-generic for $P_{\eta_{i}}$).
\end{itemize}
Since each $p_{i+1}\in M$, its initial segments are automatically in the
corresponding $N_{\xi}$'s. Then $\{ p_{i} : i < \omega \}$ generates
an $M$-generic filter $g$ for $P_{\beta}$, and for each
$\eta_{i}$, $g \cap P_{\eta_{i}} = g_{i}$. Then by the induction
hypothesis, $g \cap P_{\eta} \in
N_{o.t.(M \cap \eta)}$ for all $\eta \in M \cap \beta$, since
for all $i < \omega$ with $\eta_{i} \geq \eta$, $g \cap P_{\eta} =
g_{i} \cap P_{\eta}$.
\end{proof}

%The lemma below follows by
%letting $p_{0}$ below be a potential-$M$-generic for $P_{\beta}$.

\begin{lem}\label{fixone} Fix $\beta \leq \kappa$. If for all
$p_{0} \in P_{\beta}$ there is a $(\bar{Q}, \beta,
\beta)$-system  $(M, \bar{N}, p, q)$ with $p_{0} \in M$
and $p \leq p_{0}$ a potential-$M$-generic,
then the forcing  $P_{\beta}$ is $(\omega, \infty)$-distributive.
\end{lem}

\begin{proof} Towards a contradiction, fix the least $\beta$ for
which the lemma fails, and let $p_{0}$ be a condition in
$P_{\beta}$ forcing that
$P_{\beta}$ adds a new $\omega$-sequence of ordinals.
Let $(M, \bar{N}, p, q)$ be a $(\bar{Q}, \beta, \beta)$-system
with $p_{0} \in M$ and $p \leq p_{0}$ a potential-$M$-generic,
as given by the hypothesis of the lemma. Then there exists a
$P_{\beta}$-name
$\tau$ in $M$ such that $p_{0}$ forces that $\tau$ will be a new
$\omega$-sequence
of ordinals. Let $g \subset
P_{\beta} \cap M$ be an $M$-generic filter witnessing that $p$ is a
potential-$M$-generic. We wish to see that $q$ is below each member of
$g$. Then we will be done, as for each integer $i$ there is a member of
$g$ intersecting the antichain in $P_{\beta}$ determining the $i$th
member of $\tau$. So we will show by induction on $\rho \in (\beta +1)
\cap M$
that $q \restrict \rho \leq p' \restrict \rho$ for all $p' \in g$
(simultaneously). The cases where $\rho = 0$ or $\rho$ is a limit and
$M \cap \rho$ is
cofinal in $\rho$ are clear.

For the successor step from $\rho$ to $\rho + 1$,
$\undertilde{f}_{\rho}(M \cap \omega_{1})$ is a
$P_{\rho}$-name in $N_{o.t.(M \cap \rho)}$ for a countable
ordinal, so
by Conditions \ref{syssix} and \ref{g} of Definition \ref{sysdef},
$q\restrict \rho$ forces that $$M \cap \lambda_{\rho + 1} \in
A^{\rho}_{\undertilde{f}_{\rho}(M \cap \omega_{1})}.$$
Then  $q\restrict \rho$ forces that $p(\rho) = g(\rho)^{\frown}(\cup
g(\rho))$.
By Condition \ref{syslast} of Definition \ref{sysdef},
$q \restrict (\rho + 1) \leq p\restrict (\rho+1)$, so
$q \restrict \rho$
forces that $q(\rho)$ extends
$p(\rho)$. Now fix $p' \in g$. Since $q \restrict \rho$ forces that
$p(\rho)$ extends $p'(\rho)$, we have that $q
\restrict (\rho + 1) \leq p' \restrict (\rho + 1)$.

Lastly, for the case where $\rho$ is a limit and $M \cap \rho$ is not
cofinal
in $\rho$, note that $\rho$ has uncountable cofinality, and also that
we have shown that each $P_{\rho'}$, $\rho' < \rho$ is $(\omega,
\infty)$-distributive. By Property \ref{propone} then, densely many
conditions in
$P_{\rho}$ (and therefore $g \cap P_{\rho}$) are conditions in
some $P_{\rho'}$, $\rho' < \rho$, and so densely many conditions
in $g \cap P_{\rho}$ are in some $P_{\rho'}$ with $\rho' \in
\rho \cap M$.
\end{proof}

Given a countable $X \prec (H(\chi), \in, \leq_{\chi})$ with
$\lambda \in X$, and given $\eta < \lambda$,
%and countable elementary submodels
%$X \subset Y$ of $H(\chi)$ with
%$\lambda \in X$,
we say that $Y\prec (H(\chi), \in, \leq_{\chi})$ is a
\emph{minimal} $(\eta,\lambda)$-\emph{extension} of $X$ if the
following hold.
\begin{enumerate}
\item $X \cap \eta = Y \cap \eta$. \item $Y = Sk_{(H(\chi), \in,
\leq_{\chi})}(X \cup A)$, for some $A \subset \lambda$.
\end{enumerate}
The fact about these extensions that we will use is given in the following
lemma.

\begin{lem}\label{minex} Let $X \prec (H(\chi), \in , \leq_{\chi})$, and let
$\lambda < \gamma$ be ordinals in $X$ with $\gamma$ a regular
cardinal. Let $Y$ be a minimal $(\eta, \lambda)$-extension of $X$ for
some $\eta < \lambda$. Then $X \cap
\gamma$ is cofinal in $Y \cap \gamma$.
\end{lem}

\begin{proof} Since $\gamma$ is regular, each $f\colon [\lambda]^{<\omega}
\to \gamma$ has bounded range below $\gamma$. If $f$ is in $X$
then this bound exists in $X$.
\end{proof}

\begin{lem}\label{issys} For any set $x \in H(\chi)$ and
any ordinal $\beta \leq \kappa$
there is a $(\bar{Q}, 0, \beta)$-system
$(M, \bar{N}, 1_{P_{\beta}}, \emptyset)$ with $x \in M$.
\end{lem}

\begin{proof} Let $M_{0} \prec (H(\chi), \in,
\leq_{\chi})$ be countable with $\{ \bar{Q}, \beta, x \} \in M$.
For some $\zeta \leq \omega_{1}$ we build $\langle M_{\xi} : \xi
\leq \zeta \rangle$, $\langle \gamma_{\xi} : \xi \leq \zeta
\rangle$ and $\langle N_{\xi} : \xi < \zeta \rangle$ satisfying
the following conditions.
\begin{enumerate}
\item Each $M_{\xi}$ is a countable elementary submodel of
$(H(\chi), \in, \leq_{\chi})$. \item Each $N_{\xi}$ is a countable
elementary submodel of $(H(\chi), \in, \leq_{\chi})$. \item Each
$\gamma_{\xi}$ is the $\xi$th ordinal in $M_{\xi}$ ($0$ being the
$0$th ordinal). \item $\bar{Q} \in N_{0}$. \item For all $\xi <
\zeta $, $\langle N_{\eta} : \eta < \xi \rangle \in N_{\xi}$.
\item\label{zix} For all $\xi < \zeta$, $M \cap
V_{\lambda_{\xi}+2} \in N_{\xi}$. \item\label{sis} For all $\xi <
\zeta$ and for all $\nu \in N_{\xi} \cap \omega_{1}$, $M_{\xi+1}
\cap \lambda_{\gamma_{\xi} + 1} \in A^{\gamma_{\xi}}_{\nu}$.
\item\label{seben} For each $\xi < \zeta$, $M_{\xi + 1}$ is a
minimal $(2^{2^{\lambda_{\gamma_{\xi}}}},
\lambda_{\gamma_{\xi}+1})$-extension of $M_{\xi}$, \item If $\xi$
is a limit ordinal, then $M_{\xi} = \bigcup \{ M_{\eta} : \eta <
\xi \}$. \item\label{nonna} If there exists a $\xi < \omega_{1}$
such that $\gamma_{\xi} = \beta$, then $\zeta$ is the least such
$\xi$; otherwise $\zeta = \omega_{1}$.
\end{enumerate}

Given $M_{\xi}$, by applying Condition \ref{upthre} in the
statement of Theorem \ref{both} repeatedly, once for each $A^{\gamma_{\xi}}_{\nu}$
with $\nu = N_{\xi} \cap \omega_{1}$,
we can choose each $M_{\xi+1}$ to meet Conditions \ref{sis} and \ref{seben}.
While the model $Y$ resulting from this repeated application
may not be a minimal $(2^{2^{\lambda_{\gamma_{\xi}}}},
\lambda_{\gamma_{\xi}+1})$-extension of $M_{\xi}$,
$$Y^{*} = Sk_{(H(\chi), \in, \leq_{\chi})}(M_{\xi} \cup (Y \cap
\lambda_{\gamma_{\xi}+1}))$$
will be, as $Y \cap V_{\lambda_{\gamma_{\xi}+1}} = Y^{*} \cap
V_{\lambda_{\gamma_{\xi}+1}}$, so we can take this $Y^{*}$ as
our $M_{\xi+1}$.
The rest of the construction is straightforward.

We claim that for each $\xi \leq \zeta$ the sequence $\langle
\gamma_{\eta} : \eta \leq \xi \rangle$ lists
the first $\xi + 1$ ordinals of $M_{\xi}$ in increasing order.
This follows by induction on $\xi$. It is clear when
$\xi = 0$, and when $\xi = \xi' + 1$ it follows from
Condition \ref{seben} and the fact that
$$2^{2^{\lambda_{\gamma_{\xi'}}}} > \gamma_{\xi'}$$
(so in this case $\gamma_{\xi} = \gamma_{\xi'} + 1$).
If $\xi$ is a limit ordinal, we have from the fact that
$M_{\xi} = \bigcup\{ M_{\eta} : \eta < \xi \}$
and the induction hypothesis that $\langle \gamma_{\eta} :
\eta < \xi \rangle$ lists the first $\xi$ ordinals of
$M_{\xi}$ in increasing order. Then the definition of
$\gamma_{\xi}$ finishes the proof of the claim.

% so by Condition
%\ref{seben} the sequence of $\gamma_{\xi}$'s is increasing. Also, if
%$\eta < \gamma_{\xi}$ is an ordinal in $M_{\xi}$, then there is some
%$\xi' < \xi$ such that $\eta = \gamma_{\xi'}$. This can be seen by induction.
%It is clear for the case $\xi = 0$, and for the successor step, note that
%$\gamma_{\xi + 1} = \gamma_{\xi} + 1$. If $\bar{\xi}$ is a limit ordinal,
%fix $\eta < \gamma_{\bar{\xi}}$, and fix $\xi'< \bar{\xi}$ such that
%$\eta$ is the $\xi'$-th ordinal in $M_{\bar{\xi}}$. Fix
%$\xi''< \bar{\xi}$ greater than $\xi'$ such that $\eta \in M_{\xi''}$.
%Then $\gamma_{\xi''}$ must be greater than $\eta$, since
%$M_{\xi''}$ is a subset of $M_{\bar{\xi}}$, so the
%$\xi''$-th ordinal in $M_{\xi''}$ must be at least as big as the
%$\xi''$-th ordinal in $M_{\bar{\xi}}$ (and in fact is the same, by the
%previous remark).
%Then we can apply the induction hypothesis for $M_{\xi''}$.

Next we claim that $\gamma_{\xi} = \beta$ for some
$\xi < \omega_{1}$, and so $\zeta < \omega_{1}$. Given
this we are done, letting
$M = M_{\zeta}$ and $\bar{N} = \langle N_{\xi} : \xi < \zeta \rangle$.
All the conditions of Definition \ref{sysdef} are satisfied trivially, aside
from Condition \ref{g}. Condition \ref{g} is satisfied since for each
$\eta \in M_{\zeta} \cap \beta$, there is some $\xi < \zeta$
such that $\eta = \gamma_{\xi}$ (and so $\xi = o.t.(M \cap \eta)$).
Then $N_{\xi}$ and $M_{\xi}$
were chosen to satisfy Condition \ref{g} of Definition \ref{sysdef} by
Conditions \ref{zix} and \ref{sis} of the construction, and this
relationship was preserved for
all later $M_{\xi}$ by Condition \ref{seben}.

Assume to the contrary that $\zeta = \omega_{1}$. If $\gamma_{\omega_{1}}
= \beta$, then since every ordinal in $M_{\omega_{1}} \cap \beta$ is equal
to some $\gamma_{\xi}$, there is a limit ordinal $\xi < \omega_{1}$ such that
$$\{ \gamma_{\eta} : \eta < \xi\} = M_{\xi} \cap \beta.$$ But then
$\gamma_{\xi} = \beta$, contradicting $\zeta = \omega_{1}$. So we
may assume that $\gamma_{\omega_{1}} < \beta$.  We will show that the
cofinality of $M_{\omega_{1}} \cap \gamma_{\omega_{1}}$ is countable,
which is a contradiction since $\langle \gamma_{\xi} : \xi < \omega_{1}
\rangle$ is increasing and cofinal in it.  If $\gamma_{\omega_{1}} \leq
\lambda_{\gamma_{\xi}}$ for some $\gamma_{\xi}$ with $\gamma_{\omega_{1}}
\in M_{\xi}$, then $M \cap \gamma_{\omega_{1}} = M_{\xi} \cap
\gamma_{\omega_{1}}$, which is countable. If not, $\gamma_{\omega_{1}}
= \lambda_{\gamma_{\omega_{1}}}$. If $\gamma_{\omega_{1}}$ is singular,
let $\rho$ be the cofinality of $\gamma_{\omega_{1}}$. Then $\rho <
\lambda_{\gamma_{\xi}}$ for some $\gamma_{\xi}$ with $\{\rho,
\gamma_{\omega_{1}}\} \in M_{\xi}$, and there exists a cofinal map $f
\colon \rho \to \gamma_{\omega_{1}}$ in $M_{\xi}$.  Since $\rho \cap
M_{\omega_{1}} = \rho \cap M_{\xi}$, $f[M_{\xi} \cap \rho]$ is a countable
set cofinal in $M_{\omega_{1}} \cap \gamma_{\omega_{1}}$. The last
remaining case is that $\gamma_{\omega_{1}}$ is a regular limit cardinal.
Let $\xi$ be least with $\gamma_{\omega_{1}} \in M_{\xi}$, and fix a
cofinal sequence in $M_{\xi} \cap \gamma_{\omega_{1}}$. By Lemma
\ref{minex} and Condition \ref{seben} of the construction, this sequence
is cofinal in $M_{\omega_{1}} \cap \gamma_{\omega_{1}}$.
\end{proof}

We finish by applying the following lemma to the case $\rho = 0, \gamma =
\beta$, from Lemma \ref{issys}. By Lemma \ref{fixone}, then, we
are done.

\begin{lem}\label{alldone} Fix $\beta \leq \kappa$.
Let $(M, \bar{N}, p, q)$ be a $(\bar{Q}, \rho, \beta)$-system with $\rho \leq \gamma \leq \beta$
in $M$, and assume that $p$ is a potential-$M$-generic. Then there exists
a condition $q'\in P_{\gamma}$ such that
\begin{enumerate}
\item $q' \restrict \rho = q$,
\item $(M, \bar{N}, p, q')$ is a $(\bar{Q}, \gamma, \beta)$-system.
\end{enumerate}
\end{lem}

\begin{proof}
We first note that the lemma follows from the restricted
version of the lemma where $\gamma = \beta$. To see this, fix
$M, \bar{N}, p, q, \rho, \gamma, \beta$ as given by the
hypothesis of the unrestricted version,
and note that the restricted version then
gives a $q^{*} \in P_{\gamma}$ such that
$q^{*} \restrict \rho = q$ and
$(M, \bar{N} \restrict o.t.(M \cap \gamma), p \restrict \gamma,
q^{*})$ is a $(\bar{Q}, \gamma, \gamma)$-system.
Now, if $G_{\rho} \subset P_{\rho}$ is
$V$-generic with $q \in G_{\rho}$, then $M[G_{\rho}] \cap Ord = M \cap Ord$
(since $p$ is a potential-$M$-generic and $q \leq p$)
and the following hold.
%$\bar{N}[G_{\rho}] = \langle
%N_{\gamma}[G_{\rho}] : \gamma < o.t.(M \cap \beta) \rangle \in
%SEQ_{o.t.(M \cap \beta)}(\chi)^{V[G_{\rho}]}$.
\begin{itemize}
\item $\forall \alpha < o.t.(\beta \cap M)\text{ }N_{\alpha}[G_{\rho}]
\cap \omega_{1} = N_{\alpha} \cap \omega_{1}$,
\item $\bar{N}[G_{\rho}] = \langle
N_{\alpha}[G_{\rho}] : \alpha < o.t.(M \cap \beta) \rangle \in
SEQ_{o.t.(M \cap \beta)}(\chi)^{V[G_{\rho}]}.$
\end{itemize}
Furthermore, in $V[G_{\rho}]$,
\begin{itemize}
\item $q^{*}\restrict [\rho, \gamma) \in P_{\gamma}/P_{\rho}$
%\item $q^{*}$ is $(\bar{N}[G_{\rho}]\restrict (o.t.(M \cap \rho),
%o.t.(M \cap \gamma)), P_{\gamma}/P_{\rho})$-semi-generic),
\item $\forall\gamma' \in (\gamma\setminus \rho) \cap M$,
$q^{*} \restrict [\rho, \gamma')$
is $(N_{o.t.(M \cap \gamma')}[G_{\rho}], P_{\gamma'}/P_{\rho})$-semi-generic,
\item $p \restrict [\rho, \gamma) \geq q^{*} \restrict [\rho, \gamma)$.
\end{itemize}
Applying the elementarity of
$N_{o.t.(M \cap \gamma)}[G_{\rho}]$ in $H(\chi)^{V[G_{\rho}]}$, we see
that there exists a condition $\bar{q}\in N_{o.t.(M \cap \gamma)}[G_{\rho}]\cap
P_{\gamma}/P_{\rho}$ satisfying
these conditions.
By part \ref{tweetie} of Lemma \ref{cspres}, then,
there is a condition $q' \in P_{\gamma}$ such that $q' \restrict
\rho = q$ and $q$
forces
that $q' \restrict [\rho, \gamma)$ is a $$(\bar{N}[G_{\rho}]\restrict
[o.t.(M \cap \gamma),
o.t.(M \cap \beta)), P_{\gamma}/P_{\rho})\text{-semi-generic}$$
condition extending such a $\bar{q}$.
This $q'$ suffices.
We have Conditions \ref{syssix} and \ref{syslast}
of Definition \ref{sysdef} by the properites listed above for
$q^{*}\restrict [\rho, \gamma)$,
Condition \ref{sysfive} by our extension, and the others by the assumptions of
the lemma.

The restricted version of the lemma
follows by induction on $o.t.((\beta \setminus \rho) \cap M)$.
Note that our induction hypothesis entitles us (once
we have fixed $\beta$ and $\rho$) to
assume that the unrestricted version holds whenever
$\gamma < \beta$.

Now, when $\rho = \beta$ there is  nothing to show.
The argument for increasing $o.t.((\beta\setminus \rho)\cap M)$ by
one follows from the case where $\beta = \rho + 1$, and this
case is also trivial (see Remark \ref{sdog}).

%This requires only extending
%the $$(\bar{N}\restrict [o.t.(M \cap \rho), o.t.(M \cap \beta)),
%P_{\rho})\text{-semi-generic}$$
%$q$ to a $$(\bar{N} \restrict [o.t.(M \cap \gamma), o.t.(M \cap
%\beta)),P_{\gamma})\text{-semi-generic}$$
%$q'$ with $q' \restrict \rho = q$, this case follows from
%part \ref{tweetie} of Theorem \ref{cspres}.
%exists in $N_{o.t.(M \cap \gamma)}[G_{\rho}]$.
%That is, for each $\delta$ in $(\rho, \gamma]$, $q^{*} \restrict
%(\rho, \delta)$ is
%$(\bar{N}[G_{\rho}]\restrict [o.t.(M \cap \delta),
%o.t.(M \cap \gamma), P_{\delta}/P_{\rho})$-semi-generic.
%Using the fact that $P_{\gamma}/P_{\rho}$ is $\eta$-semi-proper
%for all countable $\eta$, extend $\bar{q}$ to a
%$$(\bar{N}[G_{\rho}]\restrict [o.t.(M \cap \gamma),
%o.t.(M \cap \beta)), P_{\gamma}/P_{\rho})\text{-semi-generic}$$
%condition $q^{**}$. Finally, we
%let $q'\restrict [\rho, \gamma)$ be a $P_{\rho}$-name below $q$ for such a $q^{**}$.
%Again letting $q' \restrict \rho = q$ as required,

The only remaining case then is when $\rho < \beta$ and $\beta =
\gamma$ is a limit ordinal. Fix an increasing sequence $\eta_{i}$
$(i < \omega)$ cofinal in $M \cap \beta$ with $\eta_{0} = \rho$,
and succesively apply the induction hypothesis (in $V$, as opposed
to some submodel). That is, let $q_{0} = q$, and given $q_{i}$,
let $q_{i+1}\in P_{\eta_{i+1}}$ be such that $q_{i+1} \restrict
\eta_{i} = q_{i}$ and $(M, \bar{N}, p, q_{i+1})$ is a $(\bar{Q},
\eta_{i+1}, \beta)$-system. Then each $q_{i}$ is
$(\bar{N}\restrict [o.t.(M \cap \eta_{i}), o.t.(M \cap \beta)),
P_{\eta_{i}})$-semi-generic, and the limit of the $q_{i}$'s is the
desired $q'$. Condition \ref{sysfive} of Definition \ref{sysdef}
is then easily satisfied (the second half being vacuous), and
Conditions \ref{syssix} and \ref{syslast}, being local properties,
are satisfied by the induction hypothesis.
\end{proof}

\section{Applications}

\subsection{Bounding}

The following is a a minor modification of a standard fact.

\begin{thrm}\label{std} Let $\lambda$ be a measurable cardinal,
let $\theta$ be a regular cardinal such that $V_{\lambda+2}
\subset H(\theta)$, and let $\leq^{*}$ be a wellordering of
$H(\theta)$. Then for every countable ordinal $\beta$ and every $X
\prec (H(\theta), \in, \leq^{*})$, there is a $z \in
[\lambda]^{\beta}$ such that
$$(**)\text{  }[z]^{<\omega} \subset \bigcap
\{ D^{Z}_{\lambda} : Z \prec (H(\theta), \in, \leq^{*}) \wedge Z
\in X\},$$ and such that, letting $$Y = Sk_{(H(\theta), \in,
\leq^{*})}(X \cup z),$$ $X \cap \lambda$ is a proper initial
segment of $Y \cap \lambda$.
\end{thrm}

\begin{proof} It suffices to prove the theorem for the case $\beta = 1$,
since then we may repeat the construction any countable ordinal
number of times. To see that $(**)$ is preserved, note that if
after repeating the construction some countable number of times we
have the $\lambda$-end-extension $X'$ of $X$ and we have added
$\gamma_{0},\ldots,\gamma_{n}$ (possibly among other ordinals) to
$z$ and $Z \prec (H(\theta), \in, \leq^{*})$ is in $X$, then $Z' =
Sk_{(H(\theta), \in, \leq^{*})}(Z \cup \{ \gamma_{0},\ldots,
\gamma_{n}\}) \in X'$. Further, if $\gamma_{n+1} \in
D^{Z'}_{\lambda}$, then by Lemma \ref{helper}, $\{
\gamma_{0},\ldots,\gamma_{n+1}\} \in D^{Z}_{\lambda}$.

Now, to prove the theorem for the case $\beta = 1$, let $\mu$ be
the $\leq^{*}$-least normal measure on $\lambda$, and fix $X$ as
in the statement of the theorem. Let $\gamma$ be any member of
$\bigcap (X \cap \mu)$, and let
$$Y = Sk_{(H(\theta), \in, \leq^{*})}(X \cup \{\gamma\}).$$
The key point is that if $f \colon \lambda \to \lambda$ is a
regressive function in any elementary submodel of $(H(\theta),
\in, \leq^{*})$, then $f$ is constant on a set $A \in \mu$, and
since $\mu$ is $\leq^{*}$-least, $\mu$ and therefore $A$ and the
constant value are also in the model. Along with the fact that
$\bigcap (X \cap \mu) \subset \bigcap (Z \cap \mu)$ for any $Z
\prec (H(\theta), \in, \leq^{*})$ with $Z \in X$, this shows that
$\gamma$ satisfies $(**)$. To see that $Y\cap \lambda$ is a proper
end-extension of $X \cap \lambda$, note that $\gamma \not\in X$,
and since any ordinal $\eta \in Y \cap \lambda$ is of the form
$f(\gamma)$ for some $f \colon \lambda \to \lambda$ in $X$, if
$\eta \in Y \cap \gamma$ then any corresponding $f$ is regressive
on a set in $\mu$, and so $\eta \in X$, by the key point.
\end{proof}

Putting together Theorems \ref{both} and \ref{std} we have the following.

\begin{cor} Say that there exists a strongly inaccessible limit of
measurable cardinals. Then there is a semi-proper forcing extension in
which Bounding holds, along with the Continuum Hypothesis.
\end{cor}

\begin{proof} Let $\langle \lambda_{\rho} : \rho < \kappa
\rangle$ be a continuous increasing sequence of cardinals with
supremum $\kappa$ strongly inaccessible such that each
$\lambda_{\rho+1}$ is a measurable cardinal. Fix a regular
cardinal $\chi > (2^{\kappa})^{+}$,  and let $\leq_{\chi}$ be a
wellordering of $H(\chi)$. Let $$\mathcal{A} = \langle
A^{\rho}_{\beta} \subset [\lambda_{\rho+1}]^{<\omega_{1}} : \rho <
\kappa, \beta < \omega_{1} \rangle$$ be such that each
$A^{\rho}_{\beta}$ is the set of countable subsets of
$\lambda_{\rho + 1}$ of ordertype greater than $\beta$. Then by
Theorem \ref{std}, Conditions \ref{up}-\ref{upthre} of Theorem
\ref{both} are satisfied (Condition \ref{upto} by the fact that
small forcing preserves measurable cardinals), so any forcing $P$
as in the statement of Theorem \ref{both} is semi-proper and
$(\omega, \infty)$-distributive. We may further require that for
every function $f$ added by some initial stage of the iteration
there is a $\rho < \kappa$ such that $f = \undertilde{f}_{\rho}$.
Each $Q_{\rho, f}$ forces that $|\lambda_{\rho + 1}| = \aleph_{1}$
and that any canonical function for $\lambda_{\rho+1}$ dominates
$f$ on a club, so $P$ forces Bounding. Further, since $2^{\omega}
< \kappa$ in $V$ and $P$ forces $\kappa = \omega_{2}$, $P$ forces
CH.
\end{proof}

Deiser and Donder have recently shown \cite{DeDo} that Bounding is equiconsistent
with a strongly inaccessible limit of measurable cardinals.

\subsection{Suslin Bounding}

The previous application can be generalized to show that the following
statement can be forced without adding reals, answering a question in
\cite{W}.

\begin{df}(\cite{FMW}) Suppose that $A \subset \mathbb{R}$. Then
$A$ is \emph{universally Baire} if for any compact Hausdorff space $X$ and any
continuous function
$\pi \colon X \to \mathbb{R}$,
the set $\{ a \in X \mid \pi(a) \in A \}$ has the property of Baire
in $X$.
\end{df}

\begin{df} \emph{Suslin Bounding} is the following statement. Suppose that
$A \subset \mathbb{R}$ is universally Baire, and fix a function
$f\colon \omega_{1} \to A$. Then there is a tree $T$ on $\omega
\times \omega_{1}$ such that $A = p[T]$ and such that $\{ \alpha <
\omega_{1} \mid f(\alpha) \in p[T \restrict (\omega \times
\alpha)]\}$ contains a club.
\end{df}

\begin{thrm} Bounding is equivalent to Suslin Bounding for
$\Pi^{1}_{1}$ sets.
\end{thrm}

\begin{proof}

For the forward direction, we adapt the proof that $\Pi^{1}_{1}$
sets are projections of trees on $\omega \times \omega_{1}$ (see
\cite{Mo}). Let $S$ be a tree on $\omega \times \omega$, let $B =
p[T]$ and let $A = \omega^{\omega} \setminus B$. Recall that for
all $x \in A$, $$S_{x} = \{ \sigma \in \omega^{<\omega} \mid (x
\restrict |\sigma|, \sigma) \in S\}$$ is wellfounded. For each $x
\in A$, let $\beta_{x}$ be the least ordinal (necessarily
countable) such that there exists a function $r_{x} \colon
\omega^{<\omega} \to \beta_{x}$ with the property that if $\sigma,
\sigma' \in S_{x}$ are such that $\sigma$ is a proper initial
segment of $\sigma'$, then $r_{x}(\sigma) > r_{x}(\sigma')$. Fix
$f \colon \omega_{1} \to A$, and for each $\alpha < \omega_{1}$
let $h(\alpha) = \beta_{f(\alpha)}$. Now fix $\gamma < \omega_{2}$
and a bijection $g \colon \omega_{1} \to \beta$ such that for
every $\alpha$ in a fixed club $C \subset \omega_{1}$,
$o.t.(g[\alpha]) > h(\alpha)$. Let $\sigma \colon \omega \to
\omega^{<\omega}$ be a bijection such that $|\sigma(i)| \leq i+1$
for all $i< \omega$. Let $T$ be the tree on $\omega \times
\omega_{1}$ consisting of all pairs $(\rho, \tau)$ such that if
$i,j < |\rho|$ $(= |\tau|)$ and $\sigma(i)$ is a proper initial
segment of $\sigma(j)$ with $$(\rho \restrict |\sigma(i)|,
\sigma(i)), (\rho \restrict |\sigma(j)|, \sigma(j))$$ both in $S$,
then $g(\tau(i)) > g(\tau(j))$. Then $p[T] = A$, and, for each $x
\in A$, if $\delta_{x}$ is the least (again, necessarily
countable) ordinal $\delta$ such that $o.t.(g[\delta]) \geq
\beta_{x}$ and $i_{x}$ is an order preserving embedding of
$\beta_{x}$ into $g[\delta_{x}]$, then $(x, g^{-1} \circ i_{x}
\circ r_{x} \circ \sigma)$ is a path through $T \restrict (\omega
\times \delta_{x})$. Since $\delta_{f(\alpha)} \leq \alpha$ for
all $\alpha \in C$, this shows that $T$ satisfies the definition
of Suslin Bounding for $f$.

For the other direction, let $W\subset \omega^{\omega}$ be the set
of functions coding wellordings of $\omega$ under some fixed
coding with the property that for each $g \in \omega^{\omega}$ and
$n \in \omega$, $g \restrict (n+1)$ codes how $n+1$ compares with
each $m < n+1$. Each function from $\omega_{1}$ to $\omega_{1}$
induces a corresponding function from $\omega_{1}$ into $W$. Let
$f \colon \omega_{1} \to W$ be such a function, and let $T$ be the
tree given by Suslin Bounding, with $C$ the witnessing club set.
We have a partial, not necessarily transitive, order $\leq_{o}$ on
the sequences in $T$, where $\sigma_{0} \leq_{o} \sigma_{1}$ if
one extends the other, and, $\sigma_{i}$ being the longer one, the
first cordinate of the last element of $\sigma_{i}$ codes that
$|\sigma_{0}| < |\sigma_{1}|$ in the corresponding ordering. Let
$\leq_{t}$ be the least transitive ordering containing $\leq_{o}$.
Then $\leq_{t}$ is a wellfounded partial order. Seeing this
requries checking that the resulting order is antireflexive and
wellfounded. To see antireflexivity, assume that $$\tau \leq_{o}
\sigma_{0} \leq_{o} \ldots \leq_{o} \sigma_{n} \leq_{o} \tau$$ is
the shortest possible counterexample. First note that $\sigma_{0}$
and $\sigma_{n}$ cannot be comparable in $\leq_{o}$ (or identical)
since then there would be a shorter counterexample, removing the
$\tau$'s and placing either $\sigma_{n}$ and the beginning of the
sequence or $\sigma_{0}$ at the end if they are unequal.
%the longest of $\tau, \sigma_{0}, \sigma_{n}$ would
%code a linear ordering containing the other two.
Therefore
$\sigma_{0}$ and $\sigma_{n}$ must be incompatible extensions of
$\tau$, and $\sigma_{1}$ must be comparable with $\tau$ and
distinct from $\sigma_{n}$. But then depending on whether
$\sigma_{1} \geq_{o} \tau$ or $\tau \geq_{o} \sigma_{1}$, there is
a shorter counterexample, removing either $\sigma_{0}$ or
$\sigma_{2}\ldots \sigma_{n}$ from the original sequence. To see
wellfoundedness, let $\sigma_{i}$ ($i < \omega$) be a descending
sequence in $\leq_{o}$. First note that if some $\tau \in T$ has
infinitely many extensions in the sequence, then all but finitely
many of them must be extensions of a fixed immediate successor of
$\tau$, since otherwise initial segments of $\tau$ are visited
infinitely often by the sequence, which is impossible (using
antireflexivity of $\leq_{t}$), there being only finitely many of
them. But the empty sequence has infinitely many extension in the
sequence, which means that we can build an infinite chain though
$T$ all of whose members have this property. We claim that
infinitely many members of this chain must be in the sequence,
which gives a contradiction since the chain codes a wellordering.
To see the claim, fix a member $\tau_{0}$ of the chain, and an
arbitrary integer $n$. Then there is some $\sigma_{i}$, $i > n$,
extending $\tau_{0}$, and a member $\tau_{1}$ of the chain with
length greater than $|\sigma_{i}|$. Let $\sigma_{j}$ be an
extension of $\tau_{1}$, for some $j > i$. Then if $\sigma_{i}$ is
not on the chain, there must be some $k$ in the interval $(i,j)$
such that $\sigma_{k}$ is an initial segment of $\tau_{1}$. Since
$n$ was arbitrary, the claim follows.

Now extend $\leq_{t}$ to a wellordering $\leq_{T}$ of $T$, and let
$\gamma$ be the length of $\leq_{T}$. Let $h \colon \omega_{1} \to
T$ be a bijection, and define $g \colon \omega_{1} \to \omega_{1}$
by letting $g(\alpha)$ be the ordertype of $\leq_{T}$ restricted
to $h[\alpha]$. Then $g$ is a canonical function for $\gamma$.
Furthermore, for a club $C' \subset C$ of $\alpha < \omega_{1}$,
$T \restrict (\omega \times \alpha) = h[\alpha]$. For these
$\alpha$, $g(\alpha)$ is greater than the ordertype of every
wellordering in the projection of $T \restrict (\omega \times
\alpha)$, and thus greater than $f(\alpha)$.
\end{proof}

\begin{remark} Note that Suslin Bounding for $\Pi^{1}_{1}$ sets and
Suslin Bounding for $\undertilde{\Pi}^{1}_{1}$ sets are identical.
\end{remark}

Instead of working directly with the definition of universal Baireness,
we will consider the equivalent (in the presence of large cardinals) form
given by Theorem \ref{eqe}. Given a
set $X$, we let $m(X)$
denote the set of countably complete ultrafilters on $X$. Recall that a
sequence of measures $\langle \mu_{i} :
i < \omega \rangle$ such that each $\mu_{i}$ concentrates on
$\kappa^{i}$ is a \emph{countably complete tower} if for any
sequence $\langle A_{i} : i < \omega \rangle$ such that each $A_{i} \in
\mu_{i}$ there is a sequence
$z \in \kappa^{\omega}$
such that each $z \restrict i \in A_{i}$.
See \cite{MSt, W} for more detail.

\begin{df} Suppose that $\kappa$ is a nonzero ordinal and that $T$ is
a tree on $\omega \times \kappa$. Then $T$ is $\delta$-\emph{homogeneous}
if there is a partial function $\pi \colon \omega^{<\omega} \to
m(\kappa^{<\omega})$ such that
\begin{enumerate}
\item if $s \in dom(\pi)$ then $\pi(s)$ is a $\delta$-complete
measure on $\kappa^{|s|}$ and
$\pi(s)(T_{s}) = 1$, where $T_{s} = \{ t \in \kappa^{|s|} :
(s,t) \in T \}$,
\item for all $x \in \omega^{\omega}$, $x \in p[T]$ if and only if
\begin{enumerate}
\item $\{ x\restrict k : k \in \omega \} \subset dom(\pi)$,
\item $\langle \pi(x \restrict k) : k \in \omega \rangle$ is a countably
complete tower.
\end{enumerate}
\end{enumerate}
\end{df}

A set $A \subset \mathbb{R}$ is $\delta$-homogeneously Suslin if
$A = p[T]$ for some $\delta$-homogeneous tree $T$. $A$ is
$^{\infty}$-homogeneously Suslin if it is $\delta$-homogeneously Suslin
for arbitrarily large $\delta$.

\begin{thrm}\label{eqe} (\cite{FMW}) Suppose that there is a proper class of
Woodin cardinals and that $A \subset \mathbb{R}$. Then the following are
equivalent.
\begin{enumerate}
\item $A$ is universally Baire.
\item $A$ is $^{\infty}$-homogeneously Suslin.
\end{enumerate}
\end{thrm}

We use the following fact to ensure that our forcing
iteration considers all universally Baire sets.

\begin{lem} Let $P$ be an $(\omega, \infty)$-distributive
partial order, and let $G \subset P$ be $V$-generic.
Then in $V[G]$, for all $A \subset \mathbb{R}$ and all
$\gamma > (2^{|P|})^{+}$, $A$
is $\gamma$-homogeneously Suslin if and only if
$A \in V$ and $A$ is $\gamma$-homogeneously Suslin in $V$.
\end{lem}

\begin{proof}
This follows from the following standard facts about measures, where
$V[G]$ is an extension by a forcing $P$ such that $(2^{|P|})^{+}
< \gamma \leq \kappa$.
\begin{enumerate}
\item For every $\gamma$-complete measure $U$ on $\kappa^{<\omega}$ in $V[G]$,
$U \cap V \in V$. (Otherwise, densely often in $P$ there is a set in
$V$ whose membership in the measure is undecided; by genericity
then there will be a subset of the measure of size $\leq |P|$ with
empty intersection.)
\item For every $\gamma$-complete measure on $\kappa^{<\omega}$ in
$V[G]$, every positive set contains a positive set in $V$. (For
each $P$-name for a positive set and each condition in $P$, consider
the set of sequences $p$ forces into the positive set.)
\item Every $\gamma$-complete measure on $\kappa^{<\omega}$ in $V$
extends to one in $V[G]$. (All sets containing positive sets from
the ground model.)
\end{enumerate}

For the forward direction, let $T$ be a $\gamma$-homogeneous
tree on $\omega \times \kappa$ in $V[G]$ such that
$p[T] = A$, for some $\gamma > (2^{|P|})^{+}$.
Let $\pi \colon\omega^{<\omega} \to
m(\kappa^{<\omega})$ witness that $T$ is $\gamma$-homogeneous. Each
$\pi(\sigma)$ extends a measure in $V$ on $\kappa^{<\omega}$, and
since $P$ is $(\omega, \infty)$-distributive, the corresponding function
$\pi'$ taking each $\sigma$ to the restriction of $\pi(\sigma)$ to $V$
exists in $V$. For each $x \in \mathbb{R} \setminus p[T]$, let
$\langle A^{x}_{k} : k < \omega \rangle$ be a witness to the fact
that $\langle \pi(x \restrict k) : k < \omega \rangle$ is not countably
complete.  For each $s \in \omega^{<\omega}$, let
$$B_{s} = \bigcap\{ A^{x}_{|s|}
\mid x \in \mathbb{R} \setminus p[T] \wedge s
\subset x \}.$$
Since every positive set for each $\pi(x \restrict k)$ contains
one from $V$, and since $P$ is $(\omega, \infty)$-distributive, we can
assume by shrinking if necessary that $\langle B_{s} \mid s \in
\omega^{<\omega} \rangle$ is in
$V$.
Now let $T'\in V$ be the set of pairs $s,t$ such that $t \in
B_{s}$.  Since the measures are all $\gamma$-complete, each $B_{s}$ is
positive for $\pi(s)$. Then $T'$ is $\gamma$-homogeneous (with
$\pi'$ as a witness) with the same projection as $T$.

For the other direction, assume that $T$ (on $\omega \times
\kappa$) and $\pi$
in $V$ witness that $A$ is $\gamma$-homogeneously Suslin.
Extend the $\pi(\sigma)$'s to $V[G]$-measures, inducing a function
$\pi' \colon \omega^{<\omega} \to m(\kappa^{<\omega})$.
Since each positive set in $V[G]$ contains one in $V$, and since
no $\omega$-sequences of ordinals have been added by $P$, for
each $x \in \omega^{\omega}$ the countable completeness of the
corresponding tower is not changed by $P$. Since no countable sets of
ordinals have been added the projection of $T$ is the same, so $T$
(along with $\pi'$) witnesses in $V[G]$ that $A$ is $\gamma$-homogeneously
Suslin.
\end{proof}

Let $T$ be a tree on $\omega \times \kappa$ and let $f \colon
\omega_{1} \to p[T]$.
Our one-step forcing $R_{f, T}$ is the set of
continuous increasing sequences  $\langle x_{\alpha} \in
[\kappa]^{<\omega_{1}} : \alpha \leq
\beta \rangle$ of countable length such that
for all $\alpha \leq \beta$, $x_{\alpha} \cap \omega_{1} \in
\omega_{1}$ and $f(x_{\alpha} \cap \omega_{1}) \in  p[T \restrict (\omega
\times x_{\alpha})]$, ordered by extension.

Given a $\delta^{+}$-homogeneous set of reals, we use the measures
witnessing homogeneity to suitably expand countable
elementary substructures of $H(\chi)$.

\begin{lem}\label{trext} Fix $\delta \geq \omega_{1}$, and let $T$ be a
$\delta^{+}$-homogeneous tree on $\omega \times \kappa$, witnessed
by $\pi \colon \omega^{<\omega} \to m(\kappa^{<\omega})$. Let
$\chi > 2^{\kappa}$ be a regular cardinal with $\leq_{\chi}$ a
wellordering of $H(\chi)$, and let $X \prec (H(\chi), \in,
\leq_{\chi})$ with $T, \pi \in X$. Then for any countable $a
\subset p[T]$ there exists a countable $z \subset \kappa$ such
that letting $Y =  Sk_{(H(\chi), \in, \leq_{\chi})}(X \cup z)$ we
have that
\begin{enumerate}
\item $X \cap \delta = Y \cap \delta$, \item $a \subset p[T
\restrict (\omega \times (Y \cap \kappa))]$. \item $[z]^{<\omega}
\subset \bigcap \{ D^{Z}_{\kappa} : Z \prec (H(\chi), \in,
\leq_{\chi}) \wedge \{ T, \pi\} \in Z \in X\}.$
\end{enumerate}
\end{lem}

\begin{proof} Fix $T, \pi, X$ and $a = \{ a_{i} : i < \omega \}$. By Lemma \ref{helper},
it suffices to show that we can deal with $a_{0}$, as we can just
repeat the process $\omega$ times. For each $k < \omega$, let
$A_{k} = \bigcap (X \cap \pi(a_{0}\restrict k))$. Then by the
definition of $\delta^{+}$-homogeneity (i.e., countable
completeness) there exists $z \in (\kappa \setminus
\delta)^{\omega}$ such that for all $k < \omega$, $(a_{0}
\restrict k, z \restrict k) \in T$  and $z \restrict k \in A_{k}$.
Now $z$ is as desired, since by the $\delta^{+}$-completeness of
each $\pi(a_{0} \restrict k)$, if $k \in \omega$ and $h\colon
\kappa^{k} \to \delta$ then $h$ is constant on a set $C_{h} \in
\pi(a_{0}\restrict k)$, and so this constant value must be an
element of any elementary submodel of $(H(\chi), \in,
\leq_{\chi})$ with $h$ and $\pi$ as members. Further, if $h$ is in
$X$, then $z \restrict k \in C_{h}$, so $h(z \restrict k)$ is the
corresponding constant value.
\end{proof}

Given $Y$ as in Lemma \ref{trext} where $f(Y \cap
\omega_{1}) \in a$, the union of any $Y$-generic for $R_{f, T}$,
adjoined by its union, is a condition.

To get the consistency of Suslin Bounding from Theorem \ref{both},
we start from a proper class of Woodin cardinals, and let $\kappa$
be a strongly inaccessible cardinal such that every
$^{\infty}$-homogeneously Suslin set of reals is
$^{\infty}$-homogeneously Suslin in $V_{\kappa}$. The assumption of the
Woodin cardinals is just to make $^{\infty}$-homogeneous Suslinity equal
to universal Baireness. Let $F \colon \omega_{1} \to \mathbb{R}$ be a
wellordering  of the reals. Let our bookkeeping for $P$ be such that
each pair $(A, f)$, where $A$ is a universally Baire set (from the
ground model) and
$f \colon \omega_{1} \to A$ is added by some initial segment of the
iteration, is associated to some $\rho$ greater than the stage
at which $f$ was added, such that $A$ is the projection of a
$|\mathcal{P}(\mathcal{P}(\lambda_{\rho}))|^{+}$-homogeneously Suslin tree $T \in V_{\kappa}$
on $\omega \times \kappa$ for some $\kappa \leq \lambda_{\rho + 1}$.
For this $\rho$, we let each  $A^{\rho}_{\beta}$
be the set of countable $ x \subset \lambda_{\rho + 1}$
such that $F(\beta) \in p[T \restrict (\omega \times (x\cap \kappa))]$,
so that $Q_{\rho, (F^{-1}\circ f)} = R_{f, T}$.
By Lemma \ref{trext},
these sets also satisfy Conditions \ref{up}-\ref{upthre} of
Theorem \ref{both}.
This scheme then gives the following corollary.

\begin{cor} Suppose that there is a proper class of Woodin cardinals
and let $\kappa$ strongly inaccessible be such that every
$^{\infty}$-homogeneously Suslin set of reals is
$^{\infty}$-homogeneously Suslin in $V_{\kappa}$.
Then there is a semi-proper $(\omega, \infty)$-distributive
forcing of size $\kappa$ in whose
extension the Continuum Hypothesis and Suslin Bounding
hold.
\end{cor}

\noindent{Department of Mathematics and Statistics\\ Miami University\\ Oxford, Ohio 45056\\
USA}\\

\noindent{larsonpb\@@muohio.edu}

\begin{tabbing}

Institute of Mathematics\hspace{.45 in}\=Department of Mathematics\\
%\hbox\{}\\
Hebrew University\>Rutgers University\\
%\hbox{}\\
91904 Jerusalem\>New Brunswick, NJ 08903\\
Israel\> USA
\end{tabbing}

\noindent{shelah\@@math.huji.ac.il}

\end{document}